\documentclass[12pt,a4paper,dvipdfmx]{amsart}
\UseRawInputEncoding

\usepackage{amsmath,amsfonts,amsthm,amssymb,amscd,extarrows}

\usepackage{graphicx}
\usepackage{here}   
\usepackage[stable]{footmisc}
\usepackage{cancel}
\usepackage{ascmac,fancybox}
\usepackage{cases}   
\usepackage{ulem}   
\usepackage[all]{xy}
\usepackage{url}
\usepackage{color}

\usepackage{tikz}
\usetikzlibrary{intersections,calc,arrows.meta,shadows}

\makeatletter
  
  \@addtoreset{equation}{section}
\makeatother

\newtheorem{Theorem}{Theorem}[section]

\newtheorem{Lemma}[Theorem]{Lemma}

\newtheorem{Remark}[Theorem]{Remark}
\newtheorem{Example}[Theorem]{Example}
\newtheorem{Definition}[Theorem]{Definition}

\newtheorem{Claim}[Theorem]{Claim}

\def\RMN#1{\uppercase\expandafter{\romannumeral#1}}

\newcommand{\spec}{\mathop{\rm Spec}\nolimits}
\newcommand{\proj}{\mathop{\rm Proj}\nolimits}

\newcommand{\oo}{\mathcal{O}}

\newcommand{\ff}{\mathcal{F}}

\newcommand{\bZ}{\mathbb{Z}}
\newcommand{\bR}{\mathbb{R}}
\newcommand{\bRo}{\mathbb{R}_{\ge 0}}
\newcommand{\bF}{\mathbb{F}}
\newcommand{\bC}{\mathbb{C}}
\newcommand{\bQ}{\mathbb{Q}}
\newcommand{\bP}{\mathbb{P}}
\newcommand{\bN}{\mathbb{N}}

\newcommand{\bma}{{\bf a}}
\newcommand{\bmb}{{\bf b}}
\newcommand{\bmc}{{\bf c}}

\newcommand{\bmo}{{\bf o}}
\newcommand{\bmf}{{\bf f}}

\newcommand{\maru}[1]{\textcircled{\footnotesize {#1}}}

\begin{document}
\title[Infinitely generated symbolic Rees rings]{
Infinitely generated symbolic Rees rings of positive characteristic}
\author{Kazuhiko Kurano}
\dedicatory{}
\date{}
\thanks{
The second author was supported by 
 JSPS KAKENHI Grant Number 19H00637. \\
}
\maketitle

%13A30, 14E99

\begin{abstract}
Let $X$ be a toric variety over a field $K$ determined by a triangle.
Let $Y$ be the blow-up at $(1,1)$ in $X$.
In this paper we give some criteria for finite generation of the Cox ring of $Y$
in the case where $Y$ has a curve $C$ such that $C^2\le 0$ and $C.E=1$ ($E$ is the exceptional divisor).
The natural surjection $\bZ^3 \rightarrow {\rm Cl}(X)$ gives the ring homomorphism $K[\bZ^3] \rightarrow K[{\rm Cl}(X)]$. 
We denote by $I$ the kernel of the composite map $K[x,y,z] \subset K[\bZ^3] \rightarrow K[{\rm Cl}(X)]$.
Then Cox(Y) coincides with the extended symbolic Rees ring $R'_s(I)$.
In the case where ${\rm Cl}(X)$ is torsion-free, this ideal $I$ is the defining ideal of
a space monomial curve.

Let $\Delta$ be the triangle (\ref{counterexample}) below.
Then $I$ is the ideal of $K[x,y,z]$ generated by $2$-minors of the $2 \times 3$-matrix
$\{ \{ x^7, y^2, z \}, \{y^{11}, z, x^{10} \} \}$.
(In this case, there exists a curve $C$ with $C^2=0$ and $C.E=1$.
This ideal $I$ is not a prime ideal.)
Applying our criteria,  we prove that
$R'_s(I)$ is Noetherian if and only if the characteristic of $K$ is $2$ or $3$.
\end{abstract}

\section{Introduction}

For pairwise coprime positive integers
$a$, $b$ and $c$, 
let ${\frak p}$ be the defining ideal of the space monomial
curve $(T^a, T^b, T^c)$ in $K^3$, where $K$ is a field.
The ideal ${\frak p}$ is generated by at most three binomials in $K[x,y,z]$
(Herzog~\cite{Her}).
The symbolic Rees rings of space monomial primes are deeply studied by many authors. 
Huneke~\cite{Hu} and Cutkosky~\cite{C} developed criteria 
for finite generation of such rings. 
In 1994,
Goto-Nishida-Watanabe~\cite{GNW} first found examples of infinitely generated
symbolic Rees rings of space monomial primes.
Recently, using toric geometry,
 Gonz\'alez-Karu~\cite{GK} found some sufficient conditions for 
the symbolic Rees rings of space monomial primes to be infinitely generated.

Cutkosky~\cite{C} found the geometric meaning of the symbolic Rees rings of space monomial primes. 
Let ${\Bbb P}(a,b,c)$ be the weighted projective surface with weight $a$, $b$, $c$.
Let $Y$ be the blow-up at a point in the open orbit of the toric variety ${\Bbb P}(a,b,c)$.
Then the Cox ring of $Y$ is isomorphic to the extended symbolic Rees ring
of the space monomial prime ${\frak p}$.
Therefore the symbolic Rees ring of the space monomial prime ${\frak p}$
is finitely generated if and only if the Cox ring of $Y$ is finitely generated,
that is, $Y$ is a Mori dream space.
A curve $C$ on $Y$ is called a negative curve if $C^2 < 0$ and
$C$ is different from the exceptional curve $E$. 
Cutkosky~\cite{C} proved that the symbolic Rees ring of the space monomial prime ${\frak p}$ is finitely generated if and only if
the following two conditions are satisfied:
\begin{itemize}
\item[(1)]
There exists a curve $C$ such that $C^2\le 0$ and $C\neq E$.\footnote{If $\sqrt{abc}\not\in \bQ$, this curve satisfies $C^2<0$, that is, $C$ is a negative curve.}
\item[(2)]
There exists a curve $F$ on $Y$ such that $C \cap F =\emptyset$.
\end{itemize}
In the case of ${\rm ch}(K) >0$, Cutkosky~\cite{C} proved that 
the symbolic Rees ring is Noetherian if there exists a negative curve.
In the case of ${\rm ch}(K)=0$, Inagawa-Kurano~\cite{K43} developed a very simple criterion for finite generation
in the case where a minimal generator of  ${\frak p}$ defines a negative curve $C$,
i.e., $C.E=1$.
Examples that have a negative curve $C$ with $C.E \ge 2$ 
are studied in 
Gonz\'alez-AnayaGonz\'alez-Karu~\cite{GAGK}, \cite{GAGK1} and Kurano-Nishida~\cite{KN}.

The existence of negative curves is a very difficult and important problem,
that is deeply related to Nagata's conjecture (Proposition~5.2 in Cutkosky-Kurano~\cite{CK}, Remark~\ref{imprem} (4)) and the rationality of Seshadri constant.
The existence of negative curves is studied in 
Gonz\'alez-AnayaGonz\'alez-Karu~\cite{GAGK2}, \cite{GAGK3}, Kurano-Matsuoka~\cite{KM} and Kurano~\cite{K42}.

In the case of ${\rm ch}(K) >0$,
we do not know any example such that $R_s({\frak p})$ is infinitely generated.

In this paper, 
we shall discuss finite generation in a slightly broader situation than that of the symbolic Rees ring of the defining ideal of a space monomial curve.
Now, we set up the situation dealt with in this paper and describe our results. 

Let $\Delta$ be a triangle with three vertices $(x_2,\overline{u}x_2)$, $(x_1,\overline{u}x_1)$, $(0,1)$
\begin{equation}\label{triangle}
\begin{tikzpicture}[xscale = 1, yscale = 1] 
% \node (a) at (0,0) {\dot};
% \node (b) at (2,1) {\dot};
% \node (c) at (3,-1.5) {\dot};
 \draw (0,0)--(1,1.5);
 \draw (3,-1.5)--(1,1.5);
 \draw (0,0)--(3,-1.5);
  \draw (0,0) node[anchor=east]{$(x_2,\overline{u}x_2)$};
    \draw (3,-1.5) node[anchor=west]{$(x_1,\overline{u}x_1)$};
    \draw (1,1.5) node[anchor=south]{$(0,1)$};   
 \draw (1,1.5) node{$\bullet$};  
  \draw (1,-0.5) node{$\bullet$};  
      \draw (1,-0.5) node[anchor=north]{$(0,0)$};   
      \draw (-2.5,0) node[anchor=east]{$\Delta=$};
 \end{tikzpicture}
\end{equation}
where $x_1$ and $x_2$ are rational numbers such that $x_2 \le 0 \le x_1$, $W:= x_1-x_2>0$.
Let $\overline{s}$, $\overline{t}$, $\overline{u}$ be the slopes of each edges, that is,
$\overline{s} =\frac{\overline{u}x_2-1}{x_2}$ and $\overline{t} =\frac{\overline{u}x_1-1}{x_1}$. We assume $-\infty \le \overline{t}\le -1 \le \overline{u}\le 0 \le \overline{s} \le \infty$.

Let $K$ be a field and $X$ be the toric variety determined by $\Delta$, that is,
$X=\proj E(\Delta)$ where 
\begin{equation}\label{Ehrhart}
E(\Delta) = \bigoplus_{n \ge 0}\left(\bigoplus_{(\alpha,\beta)\in n\Delta \cap \bZ^2}Kv^\alpha w^\beta\right)T^n \subset K[v^{\pm 1}, w^{\pm 1}, T]
\end{equation}
is the Ehrhart ring of $\Delta$.
Here $v$, $w$, $t$ are indeterminates.
Let $\pi:Y \rightarrow X$ be the blow-up of $X$ at $e=(1,1)$,
where $e$ is the point corresponging to the prime ideal $E(\Delta) \cap (v-1,w-1)K[v^{\pm 1}, w^{\pm 1}, T]$.

Let $E$ be the exceptional divisor.
Let $C$ be the proper transform of the curve of $X$ defined by $(w-1)T$ in $E(\Delta)$.
Then $C$ is linearly equivalent to $\pi^*\Delta - E$ and we have 
\[
C^2 = 2|\Delta| -1 = W-1 .
\]
Here remark that $C$ is isomorphic to $\bP^1_K$.
Let $u_2$($\ge 0$) and $u$($>0$) be integers such that $\overline{u} = -u_2/u$ and
$(u_2,u)=1$.
Let $\Delta'$ be the triangle with three vertices $(0,0)$, $(u,-u_2)$, 
$(-ux_2/(x_1-x_2), (u+u_2x_2)/(x_1-x_2))$.
\begin{equation}\label{Delta'}
\begin{tikzpicture}
% \node (a) at (0,0) {\dot};
% \node (b) at (2,1) {\dot};
% \node (c) at (3,-1.5) {\dot};
 \draw (0,0)--(1.2,1.8);
 \draw (3.6,-1.8)--(1.2,1.8);
 \draw (0,0)--(3.6,-1.8);
  \draw (0,0) node[anchor=east]{$(0,0)$};
    \draw (3.6,-1.8) node[anchor=west]{$(u,-u_2)$};
    \draw (1.2,1.8) node[anchor=south]{$(-ux_2/(x_1-x_2), (u+u_2x_2)/(x_1-x_2))$};   
\draw (0,0) node{$\bullet$};  
\draw (3.6,-1.8) node{$\bullet$};  
%      \draw (1,-0.5) node[anchor=north]{$(0,0)$};   
\draw (-2.5,0) node[anchor=east]{$\Delta' =$};
 \end{tikzpicture}
\end{equation}
Remark that  the slopes of edges of $\Delta'$ are  $\overline{s}$, $\overline{t}$, $\overline{u}$.
We denote the Weil divisor $\pi^*\Delta' - uE$ by $D$.
Then we have $C.D=0$.

For a positive integer $n$, we think $nC$ as a closed subscheme of $Y$ defined by
$\oo_Y(-nC)$.
We define the Cox ring of $Y$ by 
\[
{\rm Cox}(Y) = \bigoplus_{\overline{D} \in {\rm Cl}(Y)} H^0(Y, \oo_Y(D)).
\]
Even if ${\rm Cl}(Y)$ has a torsion, we can define a ring structure on ${\rm Cox}(Y)$
in this case.

We shall prove the following three theorems in Section~\ref{sec3}.

\begin{Theorem}\label{chfree}
Let $K$ be a field.
Let $\Delta$, $\Delta'$, $W$, $X$, $Y$, $C$, $D$, $\overline{s}$, $\overline{t}$, $\overline{u}$, $u_2$, $u$ be as above.
Assume $0<W \le 1$.
Then the following conditions are equivalent:
\begin{itemize}
\item[(A0)]
${\rm Cox}(Y)$ is finitely generated over $K$.
\item[(A1)]
There exists a curve $F$ on $Y$ such that
$F\cap C = \emptyset$.
\item[(A2)]
There exists a positive integer $m$ such that $\oo_Y(mD)|_{\ell C} \simeq \oo_{\ell C}$
for any positive integer $\ell$.
\item[(A3)]
There exists a positive integer $m$ such that $\oo_Y(mD)|_{muC} \simeq \oo_{muC}$.
\item[(A4)]
There exists a positive integer $m$ such that $\xi^m$ $(\in (F_{mu})^\times)$ is written as a product of elements of ${A_{mu}}^\times$ and ${\psi(B_{mu})}^\times$.
\end{itemize}
\end{Theorem}

We refer the reader to Section~\ref{prel} for definition of $\xi$, $F_{mu}$, $A_{mu}$ and $\psi(B_{mu})$.

If (A1) is satisfied, $F$ is numerically equivalent to $mD$ for some positive integer $m$.

For $i = 1, 2, \ldots, u$, we put
\[
m_i = ^\# \{ (\alpha, \beta) \in \Delta' \cap {\Bbb Z}^2 \mid \alpha = i \} .
\]
Note that $m_i \ge 1$ for all $i = 1, 2, \ldots, u$.
We sort the sequence $m_1$, $m_2$, \ldots, $m_u$ into ascending order
\[
m'_1 \le m'_2 \le \cdots \le m'_u .
\]
We say that $\Delta'$ satisfies the {\em EMU condition}  if
\[
m'_i \ge i
\]
for $i = 1, 2, \ldots, u$.

\begin{Theorem}\label{ch0}
Let $K$ be a field of characteristic $0$.
Let $\Delta$, $\Delta'$, $W$, $X$, $Y$, $C$, $D$, $\overline{s}$, $\overline{t}$, $\overline{u}$, $u_2$, $u$ be as above.
Assume $0<W \le 1$.
Then the following conditions are equivalent:
\begin{itemize}
\item[(B0)]
${\rm Cox}(Y)$ is finitely generated over $K$.
\item[(B1)]
$\oo_Y(D)|_{uC} \simeq \oo_{uC}$.
\item[(B2)]
$\xi$ $(\in (F_{u})^\times)$ is written as a product of elements of ${A_{u}}^\times$ and ${\psi(B_{u})}^\times$.
\item[(B3)]
$\Delta'$ satisfies the EMU condition.
\end{itemize}
\end{Theorem}

\begin{Theorem}\label{chp}
Let $K$ be a field of characteristic $p$, where $p$ is a prime number.
Let $\Delta$, $\Delta'$, $W$, $X$, $Y$, $C$, $D$, $\overline{s}$, $\overline{t}$, $\overline{u}$, $u_2$, $u$ be as above.
\begin{itemize}
\item[(1)]
If $0<W<1$, then ${\rm Cox}(Y)$ is finitely generated over $K$.
\item[(2)]
Assume $W = 1$.
Let $\sigma$ be the minimal positive integer such that three vertices of $\sigma \Delta$ are lattice points.
Then the following conditions are equivalent:
\begin{itemize}
\item[(C0)]
${\rm Cox}(Y)$ is finitely generated over $K$.
\item[(C1)]
There exists a non-negative integer $r$ such that $\oo_Y(\sigma p^rC)|_{\sigma p^rC}\simeq \oo_{\sigma p^rC}$.
\item[(C2)]
There exists a non-negative integer $r$ such that $\oo_Y(-\sigma p^rC)|_{\sigma p^rC}\simeq \oo_{\sigma p^rC}$.
\item[(C3)]
There exist a non-negative integer $r$ and a positive integer $j$ such that $j$ is not divided by $p$ and $H^0(\oo_Y(-\sigma jp^rC)|_{\sigma p^rC})\neq 0$.
\item[(C4)]
There exists a non-negative integer $r$ such that $H^0(\oo_Y(-\sigma p^rC)|_{\sigma p^rC})\neq 0$.
\end{itemize}
\end{itemize}
\end{Theorem}

Here, in the case $W=1$, $\sigma C$ is rationally equivalent to $(\sigma/u)D$.

We shall prove the following examples in Section~\ref{sec4}.

\begin{Example}\label{rei}
\begin{rm}
Let $g$ be a rational number such that $2 \le g \le 3$.
Let $\Delta$ be a triangle with three vertices $(g-3,\frac{3-g}{2})$, $(g-2,\frac{2-g}{2})$, $(0,1)$.
This triangle satisfies the condition in (\ref{triangle}).
In this example, $W=1$ is satisfied.
\begin{itemize}
\item[(1)]
Assume that $K$ is a field of characteristic $0$.
Then ${\rm Cox}(Y)$ is finitely generated over $K$
if and only if $\frac{7}{3} \le g \le \frac{8}{3}$.
\item[(2)]
Assume that $K$ is a field of characteristic $p$, where $p$ is a prime number.
\begin{itemize}
\item[(i)]
If $\frac{7}{3} \le g \le \frac{8}{3}$, then ${\rm Cox}(Y)$ is finitely generated over $K$.
\item[(ii)]
Suppose $g=\frac{13}{6}$.
Then ${\rm Cox}(Y)$ is finitely generated over $K$
if and only if $p=2$ or $3$.
\end{itemize}
\end{itemize}
\end{rm}
\end{Example}

In the case $g=13/6$ in the above example,
we know that ${\rm Cox}(Y)$ is isomorphic to the extended symbolic Rees ring
$R'_s(I)$ where $I$ is an ideal of $K[x,y,z]$ of the form
\[
I=I_2\left(
\begin{array}{lll}
x^7 & y^2 & z \\ y^{11} & z & x^{10}
\end{array}
\right).
\]
Here, the above ideal is not a prime ideal.

Remark that Sannai-Tanaka~\cite{ST} constructed examples of prime ideals $I$
such that symbolic Rees rings are not finitely generated over finite fields.

\begin{Remark}
\begin{rm}
Assume that $\Delta$ satisfies $W=1$.
If $I$ is a prime ideal, then we can prove that $I$ is the defining ideal of the space monomial curve $(T^1,T^1,T^1)$,
that is, $I = (x-y,y-z)$.

Therefore, if $W=1$, there does not exist infinitely generated $R_s(I)$ such that
$I$ is a space monomial prime ideal.

Gonz\'alez-AnayaGonz\'alez-Karu~\cite{GAGK3} found examples of triangles
such that $Y$ does not have a curve $C$ such that $C^2\le 0$ and $C\neq E$
in the case ${\rm ch}(K)=0$.
In this example, $W$ is a square of a rational number.
\end{rm}
\end{Remark}

\section{Preliminaries}\label{prel}

Let $\Delta$ be the triangle in (\ref{triangle}).
Let $s_2$, $s_3$, $t_3$, $t$, $u_2$, $u$ be non-negative integers such that
$\overline{s}=\frac{s_2}{s_3}$, $\overline{t}=-\frac{t}{t_3}$, $\overline{u}=-\frac{u_2}{u}$,
$(s_2,s_3)=(t,t_3)=(u_2,u)=1$.
(Here we put $s_2=1$ and $s_3=0$ if $\overline{s}=\infty$.
Similarly we put $t=1$ and $t_3=0$ if $\overline{t}=-\infty$.)

Put $\bma=(s_2,-s_3)$, $\bmb=(-t,-t_3)$, $\bmc=(u_2,u)$.
They are normal vectors of each edges of $\Delta$.
Let $a$, $b$, $c$ be pairwise coprime positive integers such that $a\bma + b\bmb+c\bmc=\bmo$.
Let $K$ be a field and $X$ be the toric variety determined by $\Delta$, that is,
$X=\proj E(\Delta)$ where $E(\Delta)$ is the Ehrhart ring of $\Delta$ as in (\ref{Ehrhart}).
We have the following diagram such that the horizontal sequence is exact:
\begin{equation}\label{classgroup}
\begin{array}{ccccccccc}
& & \bZ & & & & & & \\
& & \uparrow & \nwarrow {\scriptstyle (a \ b \ c)} & &  
\left( \begin{array}{c} \bma \\ \bmb \\ \bmc \end{array} \right) & & & \\
0 & \longleftarrow & {\rm Cl}(X)  & \longleftarrow & \bZ^3 & \longleftarrow & \bZ^2 & \longleftarrow & 0 \\
& & & & \uparrow & & & & \\
& & & & {\bN_0}^3
\end{array}
\end{equation}
Here ${\rm Cl}(X)$ is the divisor class group of $X$.
Take the semigroup rings of semigroups in the above diagram.
\[
\begin{array}{ccccc}
K[T^{\pm 1}] & =K[\bZ] &  & & \\
&\hphantom{\scriptstyle \phi_0} \uparrow {\scriptstyle \phi_0}& \hphantom{\scriptstyle \psi} \nwarrow {\scriptstyle \psi} & & \\ 
& K[{\rm Cl}(X)]  & \stackrel{\varphi}{\longleftarrow} & K[\bZ^3]  & \\
& & & \hphantom{\iota} \uparrow \iota & \\
& & &K[ {\bN_0}^3]&=K[x,y,z]
\end{array}
\]
The map $\psi\iota:K[x,y,z] \rightarrow K[T^{\pm 1}]$ is given by $\psi\iota(x)=T^a$,  $\psi\iota(y)=T^b$ and  $\psi\iota(z)=T^c$, that is, the kernel of $\psi\iota$ is the defining ideal of the space monominal curve $(T^a,T^b,T^c)$.

If the order of the torion subgroup of ${\rm Cl}(X)$ is $d$,
then ${\rm Cl}(X)$ is isomorphic to $\bZ \oplus \bZ/d\bZ$. Therefore $K[{\rm Cl}(X)]$ is isomorphic to $K[T^{\pm 1}, U]/(U^d-1)$.
Put $I = {\rm Ker}(\varphi\iota)$.
Then we know
\begin{equation}\label{SMP}
\mbox{$I$ is a prime ideal} \Longleftrightarrow \mbox{${\rm Cl}(X)$ is torsion-free}
\Longleftrightarrow \mbox{$\bZ \bma + \bZ \bmb + \bZ \bmc = \bZ^2$} .
\end{equation}

Here suppose $K = \bC$. 
We define $\phi_k: \bC[T^{\pm 1}, U]/(U^d-1) \rightarrow \bC[T^{\pm 1}]$ by
$\phi_k(T)=T$ and $\phi_k(U)=e^{2k\pi i/d}$ for $k \in \bZ$.
Then we have
\[
I = {\rm Ker}(\varphi\iota) = (\varphi\iota)^{-1}(0) =  (\varphi\iota)^{-1}(\cap_{k=0}^{d-1}{\rm Ker}(\phi_k))
= \cap_{k=0}^{d-1}{\rm Ker}(\phi_k\varphi\iota) .
\]
Here remark that each ${\rm Ker}(\phi_k\varphi\iota)$ is a prime ideal of $\bC[x,y,z]$
for each $k$.

\begin{Definition}\label{symbrees}
\begin{rm}
Let $A$ be a commutative Noetherian ring.
Let $J$ be an ideal of $A$ with minimal prime ideals $P_1$, $P_2$, \ldots, $P_d$.
We define the $n$th symbolic power of $J$ by
\[
J^{(n)} = \cap_{i=1}^d(J^nA_{P_i}\cap A) .
\]
We define
\[
R_s(J)=\oplus_{n \ge 0}J^{(n)}T^n \subset A[T]
\]
and call it the {\em symbolic Rees ring} of $J$.
We put
\[
R'_s(J)=R_s(J)[T^{-1}] \subset A[T^{\pm 1}]
\]
and call it the {\em extended symbolic Rees ring} of $J$.
\end{rm}
\end{Definition}

Let $\pi:Y \rightarrow X$ be the blow-up of $X$ at $e=(1,1)$,
where $e$ is the point corresponging to the prime ideal $E(\Delta) \cap (v-1,w-1)K[v^{\pm 1}, w^{\pm 1}, T]$.

Let $E$ be the exceptional divisor of $\pi$.

\begin{Remark}\label{imprem}
\begin{rm}
\begin{enumerate}
\item
Let $\overline{\Delta}$ be a triangle such that three vertices are rational points.
Then there exist $M \in {\rm GL}(2,\bZ)$, $r \in \bQ_{>0}$ and $\bmf \in \bQ^2$ such that $\Delta = rM\overline{\Delta} + \bmf$, where $\Delta$ is a triangle as in (\ref{triangle}).
For the proof, we use a method in Herzong~\cite{Her}.
We do not prove this result here since we do not need it in this paper.
\item
We put $t_1=t-t_3$ and $u_1 =u-u_2$.
Since $\frac{t}{t_3}\ge 1$ and $\frac{u_2}{u}\le 1$, $t_1$ and $u_1$ are non-negative integers.
Then we have
\[
I = I_2\left( 
\begin{array}{ccc}
x^{s_2} & y^{t_3} & z^{u_1} \\
y^{t_1} & z^{u_2} & x^{s_3}
\end{array}
\right)
= (x^{s_2+s_3}-y^{t_1}z^{u_1}, y^{t_1+t_3}-z^{u_2}x^{s_2}, z^{u_1+u_2}-x^{s_3}y^{t_3}) .
\]
We give an outline of the proof of it here.

We put $J = (x^{s_2+s_3}-y^{t_1}z^{u_1}, y^{t_1+t_3}-z^{u_2}x^{s_2}, z^{u_1+u_2}-x^{s_3}y^{t_3})$. 
We know $xyz$ is a non-zero divisor of $S/J$ since $S/J$ is a $1$-dimensional Cohen-Macaulay ring
by Hilbert-Burch theorem.
Therefore we have 
\[
K[ {\bN_0}^3]/J \hookrightarrow (K[ {\bN_0}^3]/J)[(xyz)^{-1}] = K[\bZ^3]/JK[\bZ^3] .
\]
Next we shall prove $K[\bZ^3]/JK[\bZ^3] = K[{\rm Cl}(X)]$.
Thus $J$ coincides with ${\rm Ker}(\varphi\iota)$.
\item
We know 
\[
{\rm Cox}(Y) = R'_s(I) 
\]
by (2.8) in \cite{K42}.
It is well-known that $R'_s(I)$ is Noetherian iff so is $R_s(I)$.
Therefore $Y$ is a Mori dream space if and only if $R_s(I)$ is finitely generated over $K$.\footnote{Remark that following conditions are equivalent:
(1) $R_s(I)$ is finitely generated over $K$, (2) $R_s(I)$ is Noetherian, (3) $R'_s(I)$ is finitely generated over $K$, (4) $R'_s(I)$ is Noetherian, (5) ${\rm Cox}(Y)$ is finitely generated over $K$, (6) ${\rm Cox}(Y)$ is Noetherian.}
\item
Let $q_1$, \ldots, $q_n$ be independent generic points in $\bP_\bC^2$.
Suppose that $n\ge 10$.
Nagata conjectured that, if a plane curve of degree $d$ passes through
each $q_i$ with multiplicity at least $r$, then $d> \sqrt{n}r$.
Nagata~\cite{Nagata} solved it affirmatively when $n$ is a square.

In the case where $I$ is a space monomial prime ideal,
the existence of negative curves is a very difficult and important problem,
that is deeply related to Nagata's conjecture (Proposition~5.2 in Cutkosky-Kurano~\cite{CK}) and the rationality of Seshadri constant.

Even if $I$ is not a prime ideal, our problem is also deeply related to Nagata's conjecture as follows; If $Y$ does not have a curve $C$ with $C^2\le 0$ except for $E$, Nagata's conjecture is true for $n=abcd$.

Suppose that $W$ is not a square of a rational number.
Under this condition, $Y$ does not have a curve $C$ with $C^2=0$.
Assume that the characteristic of $K$ is positive and ${\rm Cox}(Y)$ is not Noetherian.\footnote{We shall give an example such that the characteristic of $K$ is positive and ${\rm Cox}(Y)$ is not Noetherian in Example~\ref{rei} (2). 
However, in our example, $W=1$ and there exists a curve $C$ with $C^2=0$.}
Then there does not exists a curve $C$ such that $C^2\le 0$ and $C\neq E$.
Hence Nagata's conjecture is true for $n = abcd$.
\end{enumerate}
\end{rm}
\end{Remark}

In the rest of this section, let us recall a method in \cite{K43}.

\begin{equation}\label{StoT}
\begin{tikzpicture}[x=2cm,y=2cm]
 \filldraw[fill=lightgray,very thick] (-1.5,0.5)--(0,0)--(-1.5,2.5);
 \filldraw[fill=lightgray,very thick] (1.5,-0.5)--(0,0)--(1.5,3);   
\draw[->,>=stealth,semithick] (-1.7,0)--(1.7,0)node[below]{}; %x軸
 \draw[->,>=stealth,semithick] (0,-0.5)--(0,3)node[right]{}; %y軸
 \draw (0,0)node[below left]{$(0,0)$}; %原点 
    \coordinate (A) at (0.3,1) node at (A) [above right] {$\overline{s}$};
   \coordinate (P) at (0.7,0.5) node at (P) [below=0] {$S$};
   \coordinate (B) at (-0.7,0.7) node at (B) [below=0] {$T$};
   \coordinate (C) at (-0.3,1) node at (C) [below=0] {$\overline{t}$};
   \coordinate (D) at (0.7,-0.3) node at (D) [below=0] {$\overline{u}$};
%   \fill (C) circle [radius=2pt];
% \draw[black,very thick,domain=-1.3:1.3] plot(\x,{(\x)^2 + (1/2)})node[left]{};
% \draw[black,very thick,domain=-1.5:0.88] plot(\x,{1/(1-\x)})node[right]{$y=1/(1-x)$};
%  \draw[red,very thick,domain=-1.1:1.1] plot(\x,{1/2})node[right]{};
%  \draw[red,very thick,domain=-0.2:1.1] plot(\x,{(\x) + (1/4)})node[right]{};
%    \draw[red,very thick,domain=-1.1:0.2] plot(\x,{(-1)*(\x) + (1/4)})node[right]{};
%  \draw[red,very thick,domain=-3:3] plot(\x,{e^((-1)*(\x)^2)+0.05})node[above]{$y=e^{-x^2}$};
\end{tikzpicture}
\end{equation}
Let $S$ and $T$ be the cones in $\bR^2$ defined by 
\begin{align*}
S & = \bRo (u,-u_2) + \bRo (s_3,s_2),  \\
T & = \bRo (-u,u_2) + \bRo (-t_3,t) .
\end{align*}

\[
\begin{tikzpicture}
\draw[line width=1.5pt] (-2,0.5)--(2,-0.5);
\coordinate[label=below:O] (O) at (0,0);
 \foreach \P in {O} \fill[black] (\P) circle (0.06); 
  \coordinate  (A) at (-2,0.5); %点A
 \coordinate  (B) at (2,-0.5); %点B
 \coordinate  (C) at (2,2);
  \coordinate  (D) at (-2,2);
  \fill[lightgray] (A)--(B)--(C)--(D)--cycle;
   \draw (-1,0)node{$\overline{u}$};
  \draw (0,1)node{$Z$};
\end{tikzpicture}
\]
Let $Z$ be the cone $\bR(u,-u_2) + \bRo(0,1)$ as above.
Let $v$ and $w$ are indeterminates over $K$.
Put
\[
x = \frac{w-1}{v-1} .
\]
Here remark
\[
w=1-x+vx.
\]
Put 
\begin{align}
\nonumber K[Z] & = \bigoplus_{(\alpha,\beta) \in Z \cap \bZ^2} Kv^\alpha w^\beta \subset K[v^{\pm 1}, w^{\pm 1}] , \\
\nonumber F & = K[Z]\left[ x \right] \subset  K[v^{\pm 1}, w^{\pm 1}, \frac{1}{v-1}] , \\
\label{xalphan} x_{\alpha,n} & =v^\alpha w^{\lceil \alpha \overline{u} \rceil}x^n \in F
\end{align}
for $\alpha\in \bZ$ and $n \in \bZ_{\ge 0}$, 
where $\lceil \alpha \overline{u} \rceil$ is the least integer bigger than or equal to $\alpha \overline{u}$.
We refer the reader to Remark~4.3  in \cite{K43} for the product $x_{\alpha,n}x_{\alpha',n'}$.
Then by Proposition~4.1 in \cite{K43}, we have
\begin{align}
 F & = \bigoplus_{\alpha \in \bZ} \bigoplus_{n \ge 0} Kx_{\alpha,n} \nonumber \\
 x^\ell F & = \bigoplus_{\alpha \in \bZ} \bigoplus_{n \ge \ell} Kx_{\alpha,n} \label{F}
\end{align}
Put
\begin{align}\label{zalphan}
z_{\alpha,n} = & v^{\alpha}w^{\lceil (\alpha-n) \overline{u} \rceil} (x+x^2+x^3+\cdots)^n
= v^{\alpha}w^{\lceil (\alpha-n) \overline{u} \rceil} x^n (1+x+x^2+\cdots)^n \\
= & x_{\alpha,n}w^{\lceil (\alpha-n) \overline{u} \rceil - \lceil \alpha \overline{u} \rceil} (1+x+x^2+\cdots)^n 
\in F/x^{\ell}F
\nonumber
\end{align}
as in \cite{K43}.
We refer the reader to (4.15) in \cite{K43} for the relation between $z_{\alpha,n}$ and $x_{\alpha,n}$.
%We define subrings of $F$ as
%\begin{align*}
%A & = \bigoplus_{\alpha \ge 0} \bigoplus_{\tiny
%\begin{array}{c}
%n \ge 0 \\ (\alpha, \lceil \alpha \overline{u} \rceil+n) \in S 
%\end{array}
%} Kx_{\alpha,n} , \\
%\psi(B) & = 
%\bigoplus_{\alpha \in \bZ} \bigoplus_{\tiny
%\begin{array}{c}
%n \ge 0 \\ (\alpha-n, \lceil (\alpha-n) \overline{u} \rceil+n) \in T
%\end{array}
%} Kz_{\alpha,n} .
%\end{align*}
We put
\begin{align*}
F_\ell & = F/x^\ell F=\bigoplus_{\alpha \in \bZ} \bigoplus_{\ell > n \ge 0} Kx_{\alpha,n} , \\
A_\ell & = \bigoplus_{\alpha \ge 0} \bigoplus_{\tiny
\begin{array}{c}
\ell > n \ge 0 \\ (\alpha, \lceil \alpha \overline{u} \rceil+n) \in S 
\end{array}
} Kx_{\alpha,n} \ \subset \ F_\ell, \\
\psi(B_\ell) & = 
\bigoplus_{\alpha \in \bZ} \bigoplus_{\tiny
\begin{array}{c}
\ell > n \ge 0 \\ (\alpha-n, \lceil (\alpha-n) \overline{u} \rceil+n) \in T
\end{array}
} Kz_{\alpha,n} \ \subset \ F_\ell 
\end{align*}
as in (4.16) and (4.17) in \cite{K43}.
Remark that both $A_\ell$ and $\psi(B_\ell)$ are subrings of $F_\ell$.
Let $C$ be the proper transform of the curve of $X$ defined by $(w-1)T$
in $E(\Delta)$ (see (\ref{Ehrhart})).
Then we know that $\spec A_\ell$ and $\spec \psi(B_\ell)$ are affine open sets of $\ell C$
such that $\ell  C = \spec A_\ell \cup \spec \psi(B_\ell)$ and
$ \spec F_\ell = \spec A_\ell \cap \spec \psi(B_\ell)$.
Put $$\xi=(1-x)^u(1-x+vx)^{-u_2} \in F_\ell^\times.$$
Then $\xi$ is the transition function of the line bundle $\oo(D)|_{\ell C}$
as in (4.18) in \cite{K43}.

For integers satisfying $0 \le m \le \ell$, we put
\begin{align*}
F(m,\ell) & = x^mF/x^\ell F = \bigoplus_{\alpha \ge 0} \bigoplus_{\ell > n \ge m} Kx_{\alpha,n} \\
A(m,\ell) & = \bigoplus_{\alpha \ge 0} \bigoplus_{\tiny
\begin{array}{c}
\ell > n \ge m \\ (\alpha, \lceil \alpha \overline{u} \rceil+n) \in S 
\end{array}
} Kx_{\alpha,n} \subset F(m,\ell) , \\
B(m,\ell) & = 
\bigoplus_{\alpha \in \bZ} \bigoplus_{\tiny
\begin{array}{c}
\ell > n \ge m \\ (\alpha-n, \lceil (\alpha-n) \overline{u} \rceil+n) \in T
\end{array}
} Kz_{\alpha,n} \subset F(m,\ell)  .
\end{align*}
Remark that $A(m,\ell)$ and $B(m,\ell)$ are ideals of $A_\ell$ and $\psi(B_\ell)$, respectively.

\section{Proof of theorems}\label{sec3}

\noindent
{\it Proof of Theorem~\ref{chfree}.} \
The equivalence of (A0), (A1), (A2), (A3) are given in Theorem~3.1 and Theorem~3.2 in \cite{K43} in the case where $I$ is a prime ideal (see (\ref{SMP})).
We can prove the equivalence of them in the same way.

Next we shall prove the equivalence of (A3) and (A4).
We shall show the following claim:

\begin{Claim}\label{transition}
Let $m$ be a positive integer.
Then the following conditions are equivalent:
\begin{itemize}

\item[$\mbox{(A3)}_m$]
$\oo_Y(mD)|_{muC} \simeq \oo_{muC}$.
\item[$\mbox{(A4)}_m$]
$\xi^m$ $(\in (F_{mu})^\times)$ is written as a product of elements of ${A_{mu}}^\times$ and ${\psi(B_{mu})}^\times$.
\end{itemize}
\end{Claim}

Here we shall give an outline of the proof.
We know that $\spec A_{mu}$ and $\spec \psi(B_{mu})$ are affine open sets of ${mu} C$
such that ${mu} C = \spec A_{mu} \cup \spec \psi(B_{mu})$ and
$\spec F_{mu} = \spec A_{mu} \cap \spec \psi(B_{mu})$.
Remark that $\oo_Y(mD)|_{mu C}$ is a line bundle for any $m$
such that $\oo_Y(mD)|_{\spec A_{mu}}$ and $\oo_Y(mD)|_{\spec \psi(B_{mu})}$
are free.
Then $\xi^m$ is the transition function of $\oo(mD)|_{mu C}$
as in (4.18) in \cite{K43}.
Thus we know that $\mbox{(A3)}_m$ is equivalent to $\mbox{(A4)}_m$.

We have completed the proof of Theorem~\ref{chfree}.
\qed

\vspace{3mm}

\noindent
{\it Proof of Theorem~\ref{ch0}.} \
(B1) is equivalent to (B2) by Claim~\ref{transition}.

By Theorem~\ref{chfree}, (B1) implies (B0).
If the condition $\mbox{(A4)}_m$ in Claim~\ref{transition} is satisfied for some $m>0$,
then $\mbox{(A4)}_1$ holds in the case where the characteristic of the field $K$ is $0$
as in Proposition~5.9 in \cite{K43}.
Thus (B0) implies (B2) by Theorem~\ref{chfree} and Claim~\ref{transition}.
(In the case where $W<1$, we can prove that (B0) implies (B1) in the same way as in 
Theorem~1.1 in Kurano-Nishida~\cite{K42}.
However this method does not work in the case where $W=1$.)

The equivalence of (B2) and (B3) can be proved in the same way as in Theorem~1.2 in \cite{K43}.
\qed

\vspace{3mm}

\noindent
{\it Proof of Theorem~\ref{chp}.} \
One can prove (1) in the same way as Cutkosky~\cite{C}.

Now we shall prove (2).
In the rest of this paper, assume $$W=x_1-x_2=1.$$
In proving this theorem, we referred to Totaro's method \cite{Totaro} of constructing nef and non semi-ample divisors on a smooth surface over a finite field.

Remark that $\oo_Y(\sigma C)$ is a line bundle over $Y$ by the definition of $\sigma$.
It is obvious that (C1) is equivalent to (C2).

Since $\sigma$ is the width of $\sigma \Delta$, $\sigma$ is a multiple of $u$.
(Remember that the slope of the bottom edge of $\Delta$ is $-\frac{u_2}{u}$ and $(u,u_2)=1$.)
We know 
\[
(\sigma/u)D \sim \sigma C .
\]
By Theorem~\ref{chfree},  (C1) implies (C0).

In order to show that (C0) implies (C1), we shall prove that the condition (A2) in 
Theorem~\ref{chfree} implies (C1) in the case where the characteristic of the field $K$ is a prime number $p$.
Assume the condition (A2) is satisfied.
Let $r \ge 0$ and $j>0$ be integers satisfying $m=jp^r$ and $(p,j)=1$.
Then we know that $\oo_Y(jp^r D)|_{\ell C} \simeq \oo_{\ell C}$
for any positive integer $\ell$.
Therefore we have $\oo_Y(\sigma jp^rC)|_{\ell C}\simeq \oo_Y((\sigma/u)jp^r D)|_{\ell C} \simeq \oo_{\ell C}$ for any positive integer $\ell$.
Putting $\ell= \sigma p^r$ we have $\oo_Y(\sigma j p^r C)|_{\sigma p^rC} \simeq \oo_{\sigma p^rC}$.
Then 
\begin{equation}\label{*}
\mbox{the order of $\oo_Y(\sigma p^r C)|_{\sigma p^rC}$ (in ${\rm Pic}(\sigma p^rC)$) divides $j$.}
\end{equation}
On the other hand, without assuming (A2), we obtain
\begin{equation}\label{pbeki}
\mbox{the order of $\oo_Y(\sigma p^r C)|_{\sigma p^rC}$ (in ${\rm Pic}(\sigma p^rC)$) 
is a power of $p$}
\end{equation}
as follows.
Consider the sequence of the natural maps
\[
{\rm Pic}(\sigma p^r C) \longrightarrow{\rm Pic}((\sigma p^r-1) C) \longrightarrow {\rm Pic}((\sigma p^r -2)C) \longrightarrow  \cdots \longrightarrow  {\rm Pic}(C) =\bZ.
\]
The image of $\oo_Y(\sigma p^r C)|_{\sigma p^rC}$ is $\oo_Y(\sigma p^r C)|_{C}$ in 
${\rm Pic}(C)$.
It is $\oo_C$ since $C^2=0$ and $C\simeq \bP_K^1$.
By the exact sequence
\[
0 \longrightarrow \oo_Y(kC)/\oo_Y((k+1)C) \longrightarrow \oo_{(k+1)C}^\times \longrightarrow  \oo_{kC}^\times  \longrightarrow 1 ,
\]
we have an exact sequence
\[
H^1(\oo_Y(kC)/\oo_Y((k+1)C)) \longrightarrow {\rm Pic}((k+1)C)
\longrightarrow {\rm Pic}(kC) 
\]
for $k \ge 1$.
Therefore each element in the kernel of  ${\rm Pic}((k+1)C)
\rightarrow {\rm Pic}(kC)$ vanishi by $p$.
Consequently (\ref{pbeki}) holds.
By (\ref{*}) and (\ref{pbeki}) we know $\oo_Y(\sigma p^r C)|_{\sigma p^rC} \simeq \oo_{\sigma p^rC}$.

The implications (C2) $\Longrightarrow$ (C4) $\Longrightarrow$ (C3) are obvious.

In the rest of this section we shall prove  (C3) $\Longrightarrow$ (C2).

Remark that $\oo_Y(nC)/\oo_Y((n-1)C)$ is a line bundle over $C$ since
both $\oo_Y(nC)$ and $\oo_Y((n-1)C)$ are locally reflexive modules.
Then the following are satisfied.

\begin{Lemma}\label{ell}
Let $n$ be an integer.
\begin{enumerate}
\item
$\oo_Y(nC)/\oo_Y((n-1)C)$ is periodic with period $\sigma$.
\item
We have $\oo_Y(nC)/\oo_Y((n-1)C) \simeq \oo_{\bP_K^1}$ if $\sigma$ divides $n$.
\item
We have $H^0(\oo_Y(nC)/\oo_Y((n-1)C))=0$ if $\sigma$ does not divide $n$.
\end{enumerate}
\end{Lemma}

We shall prove this lemma after completing the proof of Theorem~\ref{chp}.

\vspace{3mm}

Let $m$ be an integer.
By the above lemma we have exact sequences
\[
0= H^0(\oo_Y((\sigma m-n+1) C)/\oo_Y((\sigma m-n )C)) \longrightarrow 
H^0(\oo_Y(\sigma m C)|_{n C}) \longrightarrow H^0(\oo_Y(\sigma m C)|_{(n-1)C})
\]
for $n= 2, 3, \ldots,\sigma$.
Thus we know that the natural map
\begin{equation}\label{inj}
\mbox{$H^0(\oo_Y(\sigma m C)|_{\sigma C}) \longrightarrow H^0(\oo_Y(\sigma m C)|_{C})$
is injective.}
\end{equation}

Assume that (C3) is satisfied.

First assume $r=0$.
Then we have the injection
\[
0 \neq H^0(\oo_Y(-\sigma j C)|_{\sigma C}) \hookrightarrow H^0(\oo_Y(-\sigma j C)|_{C}) = H^0(\oo_{\bP_K^1})
=K
\]
as in (\ref{inj}) (see Lemma~\ref{ell}).
Therefore the non-zero section in 
$H^0(\oo_Y(-\sigma j C)|_{\sigma C})$ does not vanish
at any point of $C$.
Thus we have $\oo_Y(-\sigma j C)|_{\sigma C}=\oo_{\sigma C}$ and
the order of $\oo_Y(-\sigma C)|_{\sigma C}$ divides $j$.
Then, by (\ref{pbeki}), we know that the order of  $\oo_Y(-\sigma C)|_{\sigma C}$ is a power of $p$.
Hence the order is one and $\oo_Y(-\sigma C)|_{\sigma C}=\oo_{\sigma C}$.

Next assume $r>0$.
We may assume that 
\begin{equation}\label{katei}
H^0(\oo_Y(-\sigma  j'p^{r'} C)|_{\sigma p^{r'} C})=0
\end{equation}
for any integer $r'$ such that $0\le r'<r$ and any positive integer $j'$ which is not divided by $p$.
Consider the map
\[
H^0(\oo_Y(-\sigma j p^{r} C)|_{\sigma p^{r} C}) \longrightarrow H^0(\oo_Y(-\sigma j p^{r} C)|_{C}) .
\]
It is the composite map of
\begin{equation}\label{r_1}
H^0(\oo_Y(-\sigma j p^{r} C)|_{\sigma p^{r_1} C}) \longrightarrow H^0(\oo_Y(-\sigma j p^{r} C)|_{\sigma p^{r_1-1}C}) 
\end{equation}
for $r_1=1, 2, \ldots, r$
and
\[
H^0(\oo_Y(-\sigma j p^{r} C)|_{\sigma C}) \longrightarrow H^0(\oo_Y(-\sigma j p^{r} C)|_{C}) .
\]
The last map is injective by (\ref{inj}).
The map (\ref{r_1}) is the composition of
\[
H^0(\oo_Y(-\sigma j p^{r} C)|_{\sigma(p-i) p^{r_1-1} C}) \longrightarrow H^0(\oo_Y(-\sigma j p^{r} C)|_{\sigma(p-i-1) p^{r_1-1}C}) 
\]
for $i=0,1,2,\ldots, p-2$.
The kernel of the above map is
\begin{align*}
& H^0\left( \frac{\oo_Y(-\sigma j p^{r} C-\sigma(p-i-1) p^{r_1-1}C)}{\oo_Y(-\sigma j p^{r} C-\sigma(p-i) p^{r_1-1}C)}
\right) \\
= & H^0\left( \frac{\oo_Y(-\sigma(jp^{r-r_1+1}-p+i+1) p^{r_1-1}C)}{\oo_Y(-\sigma(jp^{r-r_1+1}-p+i) p^{r_1-1}C)}
\right) \\
= & H^0(\oo_Y(-\sigma(jp^{r-r_1+1}-p+i+1) p^{r_1-1}C)|_{\sigma p^{r_1-1}C})  =0
\end{align*}
since $r_1-1<r$ (see (\ref{katei})).
Therefore (\ref{r_1}) is injective for $r_1=1, 2, \ldots, r$.
Thus we obtain the injection
\[
0 \neq H^0(\oo_Y(-\sigma j p^{r} C)|_{\sigma p^rC}) \hookrightarrow H^0(\oo_Y(-\sigma j  p^{r}C)|_{C}) = H^0(\oo_{\bP_K^1})
=K .
\]
Therefore any non-zero section in $H^0(\oo_Y(-\sigma j p^rC)|_{\sigma p^rC})$ does not vanish
at any point of $C$.
Thus we have $\oo_Y(-\sigma j p^r C)|_{\sigma p^rC}=\oo_{\sigma p^rC}$ and
the order of $\oo_Y(-\sigma p^r C)|_{\sigma p^rC}$ divides $j$.
Then, by (\ref{pbeki}), we know the order is one and $\oo_Y(-\sigma p^rC)|_{\sigma p^rC}=\oo_{\sigma p^rC}$.

We have completed the proof of Theorem~\ref{chp}.
\qed

\vspace{2mm}

\noindent
{\it Proof of Lemma~\ref{ell}.}
Remember the cones $S$ and $T$ in (\ref{StoT}).
We define 
\begin{equation}\label{ab}
\begin{split}
a_i & = ^\# \{ (\alpha, \beta) \in S \cap {\Bbb Z}^2 \mid \alpha = i \}  \\
b_i & = ^\# \{ (\alpha, \beta) \in T \cap {\Bbb Z}^2 \mid \alpha = i \} .
\end{split}
\end{equation}
By definition, we have
\[
\cdots \ge b_{-3} \ge b_{-2} \ge b_{-1} \ge b_{0}>0<a_0 \le a_1\le a_2\le a_3 \le \cdots .
\]
Remark that $a_0=a_1=\cdots=\infty$ if $\overline{s}=\infty$,
and $b_{0}=b_{-1}=\cdots=\infty$ if $\overline{t}=-\infty$.

We put
\begin{align*}
P_A & := \left\{ (\alpha,n) \in \bZ^2 \ \left| \  
\begin{array}{l}
\alpha \ge 0, \ n \ge 0, \\
(\alpha, \lceil \alpha \overline{u} \rceil+n) \in S 
\end{array}
 \right.  \right\} 
 =
 \left\{ (\alpha,n) \in \bZ^2 \ \left| \  
\begin{array}{l}
\alpha \ge 0, \ n \ge 0, \\
a_\alpha \ge n+1 
\end{array}
 \right.  \right\} 
, \\
P_B & := 
\left\{ (\alpha,n) \in \bZ^2 \ \left| \  
\begin{array}{l}
\alpha \in \bZ, \ n \ge 0, \\
(\alpha-n, \lceil (\alpha-n) \overline{u} \rceil+n) \in T 
\end{array}
 \right.  \right\} 
 =
 \left\{ (\alpha,n) \in \bZ^2 \ \left| \  
\begin{array}{l}
\alpha \in \bZ, \ n \ge 0, \\
b_{\alpha-n} \ge n+1 
\end{array}
 \right.  \right\} .
\end{align*}

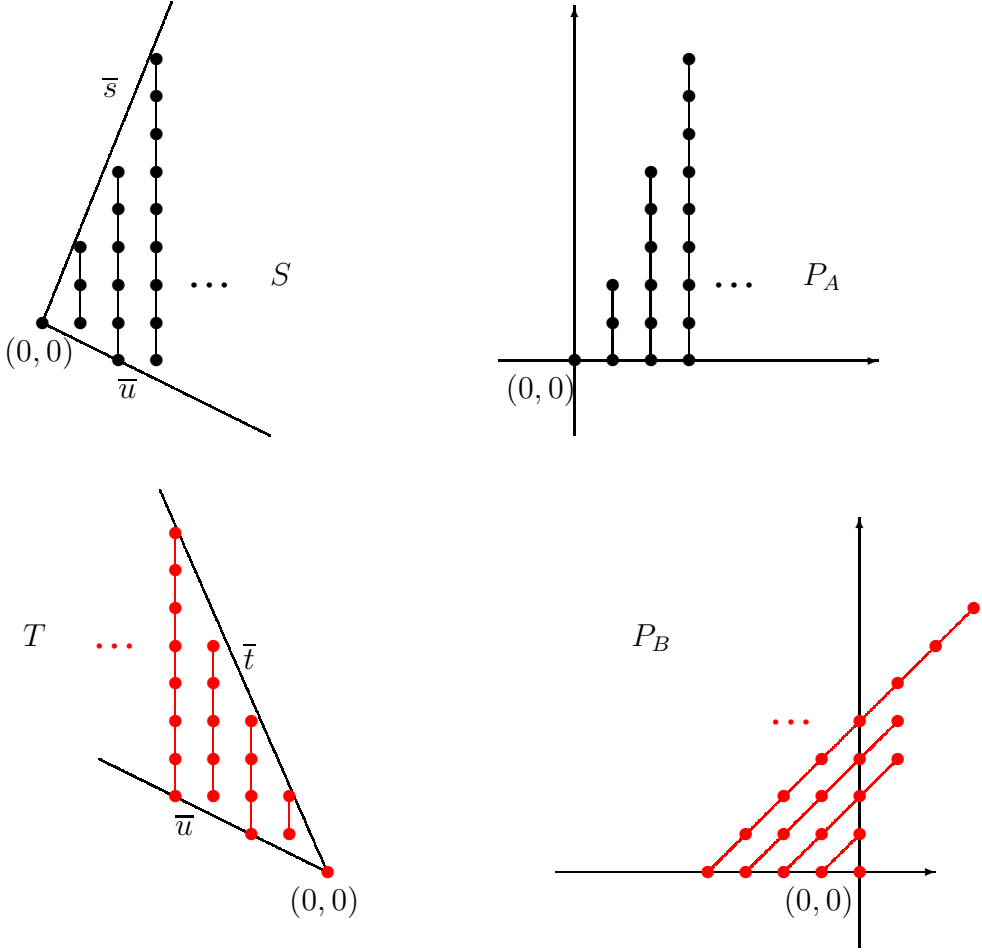
\[
{
\setlength\unitlength{1truecm}
  \begin{picture}(5.5,6)(2.5,-3)
  %\put(-1,0){\vector(1,0){7}}
  %\put(0,-3){\vector(0,1){6}}
\qbezier  (3,-1.5) (4,-2) (6,-3)
\qbezier (3,-1.5) (4,1) (4.7, 2.75)
\put(2.9,-1.6){$\bullet$}
\put(3.4,-1.6){$\bullet$}
\put(3.4,-1.1){$\bullet$}
\put(3.4,-0.6){$\bullet$}
\qbezier (3.5,-1.5) (3.5,-1) (3.5,-0.5)
\put(3.9,-2.1){$\bullet$}
\put(3.9,-1.6){$\bullet$}
\put(3.9,-1.1){$\bullet$}
\put(3.9,-0.6){$\bullet$}
\put(3.9,-0.1){$\bullet$}
\put(3.9,0.4){$\bullet$}
\qbezier  (4,-2) (4,-1.5) (4,0.5)
\put(4.4,-2.1){$\bullet$}
\put(4.4,-1.6){$\bullet$}
\put(4.4,-1.1){$\bullet$}
\put(4.4,-0.6){$\bullet$}
\put(4.4,-0.1){$\bullet$}
\put(4.4,0.4){$\bullet$}
\put(4.4,0.9){$\bullet$}
\put(4.4,1.4){$\bullet$}
\put(4.4,1.9){$\bullet$}
\qbezier (4.5,-2) (4.5,-1.5) (4.5,2)
\put(2.5,-2){$(0,0)$}
\put(6,-1){$S$}
  \put(4,-2.5){$\overline{u}$}
  \put(3.8,1.5){$\overline{s}$}
   \multiput(5,-1)(0.2,0){3}{\circle*{0.08}}
  \end{picture}
}
{
\setlength\unitlength{1truecm}
  \begin{picture}(7,6)(-2,-3)
  \put(-1,-2){\vector(1,0){5}}
  \put(0,-3){\vector(0,1){5.7}}
\put(-0.1,-2.1){$\bullet$}
\put(0.4,-2.1){$\bullet$}
\put(0.4,-1.6){$\bullet$}
\put(0.4,-1.1){$\bullet$}
\qbezier (0.5,-2) (0.5,-1.5) (0.5,-1)
\put(0.9,-2.1){$\bullet$}
\put(0.9,-1.6){$\bullet$}
\put(0.9,-1.1){$\bullet$}
\put(0.9,-0.6){$\bullet$}
\put(0.9,-0.1){$\bullet$}
\put(0.9,0.4){$\bullet$}
\qbezier (1,-2) (1,-1.5) (1,0.5)
\put(1.4,-2.1){$\bullet$}
\put(1.4,-1.6){$\bullet$}
\put(1.4,-1.1){$\bullet$}
\put(1.4,-0.6){$\bullet$}
\put(1.4,-0.1){$\bullet$}
\put(1.4,0.4){$\bullet$}
\put(1.4,0.9){$\bullet$}
\put(1.4,1.4){$\bullet$}
\put(1.4,1.9){$\bullet$}
\put(-0.9,-2.5){$(0,0)$}
\qbezier (1.5,-2) (1.5,-1.5) (1.5,2)
   \multiput(1.9,-1)(0.2,0){3}{\circle*{0.08}}
   \put(3,-1){$P_A$}
  \end{picture}
}
\]

\[
{
\setlength\unitlength{1truecm}
  \begin{picture}(7,6)(-2,-2.5)
  %\put(-1,0){\vector(1,0){7}}
  %\put(0,-3){\vector(0,1){6}}
\qbezier (0,0) (1,-0.5) (3,-1.5) 
\qbezier (3,-1.5) (2,0.8) (0.8,3.56)
  \textcolor{red}{
\put(2.9,-1.6){$\bullet$}
\put(2.4,-1.1){$\bullet$}
\put(2.4,-0.6){$\bullet$}
\qbezier (2.5,-1) (2.5,-0.7) (2.5,-0.5)
\put(1.9,-1.1){$\bullet$}
\put(1.9,-0.6){$\bullet$}
\put(1.9,-0.1){$\bullet$}
\put(1.9,0.4){$\bullet$}
\qbezier (2,-1) (2,0) (2,0.5)
\put(1.4,-0.6){$\bullet$}
\put(1.4,-0.1){$\bullet$}
\put(1.4,0.4){$\bullet$}
\put(1.4,0.9){$\bullet$}
\put(1.4,1.4){$\bullet$}
\qbezier (1.5,-0.5) (1.5,0) (1.5,1.5)
\put(0.9,-0.6){$\bullet$}
\put(0.9,-0.1){$\bullet$}
\put(0.9,0.4){$\bullet$}
\put(0.9,0.9){$\bullet$}
\put(0.9,1.4){$\bullet$}
\put(0.9,1.9){$\bullet$}
\put(0.9,2.4){$\bullet$}
\put(0.9,2.9){$\bullet$}
\qbezier (1,-0.5) (1,0) (1,3)
         \multiput(0,1.5)(0.2,0){3}{\circle*{0.08}}
}
\put(2.5,-2){$(0,0)$}
  \put(1,-1){$\overline{u}$}
  \put(1.9,1.2){$\overline{t}$}
      \put(-1,1.5){$T$}
  \end{picture}
}
{
\setlength\unitlength{1truecm}
  \begin{picture}(7,6)(-2,-3)
  \put(-1,-2){\vector(1,0){5}}
  \put(3,-3){\vector(0,1){5.7}}
  \textcolor{red}{
\put(2.9,-2.1){$\bullet$}
\put(2.4,-2.1){$\bullet$}
\put(2.9,-1.6){$\bullet$}
\qbezier (2.5,-2) (2.75,-1.75) (3,-1.5)
\put(1.9,-2.1){$\bullet$}
\put(2.4,-1.6){$\bullet$}
\put(2.9,-1.1){$\bullet$}
\put(3.4,-0.6){$\bullet$}
\qbezier (2,-2) (2.5,-1.5) (3.5,-0.5)
\put(1.4,-2.1){$\bullet$}
\put(1.9,-1.6){$\bullet$}
\put(2.4,-1.1){$\bullet$}
\put(2.9,-0.6){$\bullet$}
\put(3.4,-0.1){$\bullet$}
\qbezier (1.5,-2) (2,-1.5) (3.5,0)
\put(0.9,-2.1){$\bullet$}
\put(1.4,-1.6){$\bullet$}
\put(1.9,-1.1){$\bullet$}
\put(2.4,-0.6){$\bullet$}
\put(2.9,-0.1){$\bullet$}
\put(3.4,0.4){$\bullet$}
\put(3.9,0.9){$\bullet$}
\put(4.4,1.4){$\bullet$}
\qbezier (1,-2) (1.5,-1.5) (4.5,1.5)
 \multiput(1.9,0)(0.2,0){3}{\circle*{0.08}}
 }
 \put(2,-2.5){$(0,0)$}
   \put(0,1){$P_B$}
  \end{picture}
}
\]

Put 
\begin{equation}\label{theta}
\mbox{$\theta = \frac{-x_2\sigma}{x_1-x_2}=-x_2\sigma$ and $\theta'=\frac{x_1\sigma}{x_1-x_2}=x_1\sigma$.}
\end{equation}
Then $\theta$ and $\theta'$ are non-negative integers such that 
$\sigma=\theta+\theta'$.

\begin{Claim}\label{claim1}
\begin{enumerate}
\item
$P_A\cap P_B = \{ m(\theta, \sigma)\mid m\in \bN_0 \}$,
where $\bN_0$ denotes the set of non-negative integers.
\item
For $\alpha \in \bZ$ and $n \in \bN_0$,
$(\alpha,n)\in P_A$ if and only if $(\alpha+\theta,n+\sigma) \in P_A$.
\item
For $\alpha \in \bZ$ and $n \in \bN_0$,
$(\alpha,n)\in P_B$ if and only if $(\alpha+\theta,n+\sigma) \in P_B$.
\end{enumerate}
\end{Claim}

First we shall prove (2).
If $\overline{s}=\infty$, then $\theta=0$ and $a_0 = a_1 = \cdots = \infty$.
In this case, $P_A = \{ (\alpha,n)\mid \alpha, n\in \bN_0 \}$. 
The assertion is obvious in this case.
Assume $\overline{s}<\infty$.
Remark that $(\theta, \overline{u}\theta)$ and $(\theta, \overline{s}\theta)$ are lattice points and $\overline{s}\theta=\overline{u}\theta+\sigma$.
Therefore, for any $i \in \bN_0$, we have 
\begin{equation}\label{a}
a_{i+\theta} = a_i+\sigma .
\end{equation}
For $\alpha \in \bN_0$ and $n \in \bN_0$,
\[
(\alpha,n)\in P_A \Leftrightarrow
a_\alpha \ge n+1 \Leftrightarrow
a_{\alpha+\theta} \ge n+1+\sigma \Leftrightarrow
(\alpha+\theta,n+\sigma)\in P_A .
\]
For $\alpha <0$ and $n \in \bN_0$, then $(\alpha,n)$ and $(\alpha+\theta,n+\sigma)$ are not in $P_A$.

Next we prove (3).
Put
\[
P'_B = 
\left\{ (\alpha,n) \in \bZ^2 \ \left| \  
\begin{array}{l}
\alpha \le 0, \ n \ge 0, \\
(\alpha, \lceil \alpha \overline{u} \rceil+n) \in T 
\end{array}
 \right.  \right\} 
 .
\]
By the same way as in (2), we can prove that, for $\alpha \in \bZ$ and $n \in \bN_0$, 
$(\alpha,n)\in P'_B$ if and only if 
$(\alpha-\theta',n+\sigma)\in P'_B$.
Then, for $\alpha \in \bZ$ and $n \in \bN_0$, 
\[
(\alpha,n)\in P_B \Leftrightarrow
(\alpha-n,n)\in P'_B \Leftrightarrow
(\alpha-n-\theta',n+\sigma)\in P'_B \Leftrightarrow
(\alpha+\theta,n+\sigma)\in P_B .
\]

Now we start to prove (1).
Suppose $(\alpha,n)\in P_A\cap P_B$.
By (2) and (3), we may assume $0 \le n < \sigma$.
Since $(\alpha,n)\in P_A$, we know $\alpha\ge 0$, $a_\alpha \ge n+1$.
Since $(\alpha,n)\in P_B$, we have $b_{\alpha-n}\ge n+1$.
Therefore we know $\alpha-n\le 0$.
Put $y_1=n-\alpha\ge 0$ and $y_2=-\alpha\le 0$.
Consider the triangle $\overline{\Delta}$ such that
the slopes of edges are $\overline{s}$, $\overline{t}$, $\overline{u}$, respectively,
and $(y_1,\overline{u}y_1)$, $(y_2,\overline{u}y_2)$ are the endpoints
of the bottom edge as below.
\[
\begin{tikzpicture}[xscale = 1, yscale = 1] 
% \node (a) at (0,0) {\dot};
% \node (b) at (2,1) {\dot};
% \node (c) at (3,-1.5) {\dot};
 \draw (0,0)--(1,1.5);
 \draw (3,-1.5)--(1,1.5);
 \draw (0,0)--(3,-1.5);
  \draw (0,0) node[anchor=east]{$(y_2,\overline{u}y_2)$};
    \draw (3,-1.5) node[anchor=west]{$(y_1,\overline{u}y_1)$};
    \draw (0.3,0.5) node[anchor=south]{$\overline{s}$};   
\draw (1.8,0.5) node[anchor=south]{$\overline{t}$};   
\draw (1.1,-1.2) node[anchor=south]{$\overline{u}$}; 
 %\draw (1,1.5) node{$\bullet$};  
  \draw (1,-0.5) node{$\bullet$};  
      \draw (1,-0.5) node[anchor=south]{$(0,0)$};   
      \draw (-2.5,0) node[anchor=east]{$\overline{\Delta}=$};
 \end{tikzpicture}
\]
Since $a_\alpha\ge n+1$ and $b_{\alpha-n}\ge n+1$,
the point $(0,n)$ is contained in $\overline{\Delta}$.
Since the area of $\overline{\Delta}$ is $n^2/2$,
the point $(0,n)$ is a vertex of  $\overline{\Delta}$.
Then we know that $(y_1,\overline{u}y_1)$ and $(y_2,\overline{u}y_2)$ are lattice points.
Then $\sigma$ divides $n$.
Thus we obtain $n=\alpha=0$.
We have completed the proof of Claim~\ref{claim1}.

\vspace{3mm}

Let's go back to the proof of Lemma~\ref{ell}.
Assume $n=-q<0$.
We have $(q+1)C=\spec A_{q+1} \cup \spec \psi(B_{q+1})$ and
$\spec F_{q+1} = \spec A_{q+1} \cap \spec \psi(B_{q+1})$.
Then we have
\begin{align*}
\left. \frac{\oo_Y(-qC)}{\oo_Y(-(q+1)C)}\right|_{\spec F_{q+1} }
& = F(q,q+1) = \bigoplus_{\alpha \in \bZ} Kx_{\alpha,q} , \\
\left. \frac{\oo_Y(-qC)}{\oo_Y(-(q+1)C)}\right|_{\spec A_{q+1} }
& = A(q,q+1) = \bigoplus_{(\alpha,q)\in P_A}Kx_{\alpha,q} \subset F(q,q+1), \\
\left. \frac{\oo_Y(-qC)}{\oo_Y(-(q+1)C)}\right|_{\spec \psi(B_{q+1}) }
& = B(q,q+1) = \bigoplus_{(\alpha,q)\in P_B}Kx_{\alpha,q} \subset F(q,q+1) .
\end{align*}
Here remark that $x_{\alpha,q}=z_{\alpha,q}$ in $F(q,q+1)$ by (4.15) in \cite{K43}.
Thus we know
\begin{align*}
H^0\left( \frac{\oo_Y(-qC)}{\oo_Y(-(q+1)C)} \right)
& =
\left( \bigoplus_{(\alpha,q)\in P_A}Kx_{\alpha,q} \right) \cap \left( \bigoplus_{(\alpha,q)\in P_B}Kx_{\alpha,q} \right) , \\
H^1\left( \frac{\oo_Y(-qC)}{\oo_Y(-(q+1)C)} \right)
& =\frac{ \bigoplus_{\alpha \in \bZ} Kx_{\alpha,q} }{
 \left( \bigoplus_{ (\alpha,q)\in P_A}Kx_{\alpha,q} \right) + \left( \bigoplus_{(\alpha,q)\in P_B}Kx_{\alpha,q} \right)
} .
\end{align*}
Therefore we know 
\[
H^0\left( \frac{\oo_Y(-qC)}{\oo_Y(-(q+1)C)} \right) =
\left\{
\begin{array}{cl}
K & \mbox{(if $\sigma | q$)} ,\\
0 & \mbox{(otherwise)}
\end{array}
\right.
\]
by Claim~\ref{claim1}.
We have proved (2) and (3) in Lemma~\ref{ell} in the case $n<0$.

Recall that $\oo_Y(-\sigma C)$ is a line bundle over $Y$
such that $\oo_Y(-\sigma C)|_C \simeq \oo_{\bP_k^1}$.
Therefore, for any $n\in \bZ$, we have
\[
\oo_Y((n-\sigma)C) = \oo_Y(nC)\otimes_{\oo_Y}\oo_Y(-\sigma C)
\]
and
\begin{align*}
 & \frac{\oo_Y((n-\sigma)C)}{\oo_Y((n-\sigma-1)C)} 
= \frac{\oo_Y(nC)}{\oo_Y((n-1)C)} \otimes_{\oo_Y} \oo_Y(-\sigma C)\\
 =  & \frac{\oo_Y(nC)}{\oo_Y((n-1)C)} \otimes_{\oo_C} \oo_Y(-\sigma C)|_C
\simeq \frac{\oo_Y(nC)}{\oo_Y((n-1)C)} 
\end{align*}
by (2).
We have completed the proof of Lemma~\ref{ell}.
\qed

\section{Examples}\label{sec4}

We shall prove (1) and (2) in Example~\ref{rei} in this section.

Let $\Delta_g$ be the triangle with three vertices $(1,0)$, $(0,0)$, $(g,4)$.
By the affine transformation $\frac{1}{2}\left(\begin{array}{cc} -2 & 1 \\ 1 & 0 \end{array} \right)\left(\begin{array}{c} x \\ y \end{array} \right)+ \left(\begin{array}{c} g-2 \\ \frac{2-g}{2} \end{array} \right)$, $\Delta_g$ is transformed to the triangle $\Delta$ with three vertices $(g-3,\frac{3-g}{2})$, $(g-2,\frac{2-g}{2})$, $(0,1)$.
Remark that if $2\le g \le 3$, then $\Delta$ satisfies the conditions in (\ref{triangle}).
In this case the triangle $\Delta'$ in (\ref{Delta'}) has three vertices
$(0,0)$, $(2,-1)$, $(6-2g,g-1)$.

First we shall prove (1) in Example~\ref{rei}.
We apply Theorem~\ref{ch0} here.
Checking the EMU condition, we know that ${\rm Cox}(Y)$ is Noetherian
if and only if the point $(1,1)$ is in $\Delta'$.
By an easy calculation we know that it is equivalent to $7/3 \le g \le 8/3$.

Next we shall prove (2) in Example~\ref{rei}. 

Fix a triangle satisfying (\ref{triangle}).
We have ${\rm Cl}(Y) = {\rm Cl}(X) \oplus \bZ$, and it is independent of the base field $K$ (see (\ref{classgroup})).
For a Weil divisor $F$ on $Y$, 
let $h^0(F)_K$ denotes the dimension of $H^0(Y,\oo_Y(F))$
in the case where the base field is $K$.
It is easy to see that $h^0(F)_K$ depends only on the characteristic of $K$.
This fact implies that finite generation depends only on the characteristic of $K$.
It is easy to check $h^0(F)_{\bF_p} \ge h^0(F)_\bQ$ for any prime number $p$.
Using this inequality, we obtain that, if ${\rm Cox}(Y)$ is finitely generated in the case where the base field is of characteristic $0$, so is for any base field $K$.
The assertion (2) (i) follows from this.

In the rest of this paper, we shall prove (2) (ii) in Example~\ref{rei}. 
Suppose $g=\frac{13}{6}$.
Then $\Delta$ is the triangle with three vertices  $(-\frac{5}{6},\frac{5}{12})$,  $(\frac{1}{6},-\frac{1}{12})$, $(0,1)$.
\begin{equation}\label{counterexample}
\begin{tikzpicture}[xscale = 3, yscale = 3] 
% \node (a) at (0,0) {\dot};
% \node (b) at (2,1) {\dot};
% \node (c) at (3,-1.5) {\dot};
 \draw (-5/6, 5/12)--(1/6, -1/12);
 \draw (0,1)--(-5/6,5/12);
 \draw (1/6,-1/12)--(0,1);
  \draw (0,0) node[anchor=north east]{$(0,0)$};
    \draw (1/6,-1/12) node[anchor=west]{$(\frac{1}{6},-\frac{1}{12})$};   
        \draw (-5/6, 5/12) node[anchor=east]{$(-\frac{5}{6},\frac{5}{12})$};   
            \draw (0,1) node[anchor=south]{$(0,1)$};   
 \draw (0,0) node{$\bullet$};  
  \draw (0,1) node{$\bullet$};  
      \draw (-1.5,0.3) node[anchor=east]{$\Delta=$};
 \end{tikzpicture}
\end{equation}
In this case, $x_1=\frac{1}{6}$, $x_2 = \frac{-5}{6}$, $u=2$, $u_2=1$,
$\overline{t}=\frac{-13}{2}$, $\overline{u}=\frac{-1}{2}$, $\overline{s}=\frac{7}{10}$ in (\ref{triangle}).
In this case, $(a,b,c)=(1,1,6)$ and $d=24$.
We know $\sigma = 12$ in Theorem~\ref{chp},
and $\theta = 10$ and $\theta'=2$ in (\ref{theta}).
The sequences in (\ref{ab}) are
\[
a_0=1, \ a_1=1, \ a_2=3, \ a_3=4, \ a_4=5, \ a_5=6, \ a_6=8, \ a_7=8, \ a_8=10, \ a_9=11, \ a_{10}=13, \ \ldots 
\]
and
\[
b_0=1, \ b_{-1}=6, \ b_{-2}=13, \ \ldots .
\]
Here we have $a_{i+10}=a_i+12$ for $i\ge 0$ as in (\ref{a}) and $b_{-i-2}=b_{-i}+12$ for $i\ge 0$.
The sets $P_A$ and $P_B$ are the following:
\begin{equation}\label{pic}
\begin{tikzpicture}[xscale = 0.45, yscale = 0.45] 
\draw[color=gray] (-6.5,0) grid(21.5,26.5);
\draw[line width=0.8pt] (-6,0) grid[step=5](21.5,26.5);
\draw[->,>=stealth,semithick] (-6.5,0)--(22,0)node[below]{}; %x軸
 \draw[->,>=stealth,semithick] (0,-2)--(0,27)node[right]{}; %y軸
 \draw (0,0)node[below left]{$(0,0)$}; %原点 
 \fill (0,0) circle[radius=0.3];
  \fill (10,12) circle[radius=0.3];
   \fill (20,24) circle[radius=0.3];
 \draw (0,0) node{$\bullet$};  
 \draw (1,0) node{$\bullet$}; 
\foreach \x in {0,1,2} \draw (2,\x) node{$\bullet$};
 \draw[line width=1pt] (2,0)--(2,2);
\foreach \x in {0,1,...,3} \draw (3,\x) node{$\bullet$};
 \draw[line width=1pt]  (3,0)--(3,3);
\foreach \x in {0,1,...,4} \draw (4,\x) node{$\bullet$};
 \draw[line width=1pt]  (4,0)--(4,4);
\foreach \x in {0,1,2,3,4,5} \draw (5,\x) node{$\bullet$};
 \draw[line width=1pt]  (5,0)--(5,5);
\foreach \x in {0,1,2,3,4,5,6,7} \draw (6,\x) node{$\bullet$};
 \draw[line width=1pt]  (6,0)--(6,7);
\foreach \x in {0,1,2,3,4,5,6,7} \draw (7,\x) node{$\bullet$};
 \draw[line width=1pt]  (7,0)--(7,7);
\foreach \x in {0,1,...,9} \draw (8,\x) node{$\bullet$};
 \draw[line width=1pt]  (8,0)--(8,9);
\foreach \x in {0,1,...,10} \draw (9,\x) node{$\bullet$};
 \draw[line width=1pt]  (9,0)--(9,10);
\foreach \x in {0,1,...,12} \draw (10,\x) node{$\bullet$};
 \draw[line width=1pt]  (10,0)--(10,12);
\foreach \x in {0,1,...,12} \draw (11,\x) node{$\bullet$};
 \draw[line width=1pt]  (11,0)--(11,12);
\foreach \x in {0,1,...,14} \draw (12,\x) node{$\bullet$};
 \draw[line width=1pt]  (12,0)--(12,14);
\foreach \x in {0,1,...,15} \draw (13,\x) node{$\bullet$};
 \draw[line width=1pt]  (13,0)--(13,15);
\foreach \x in {0,1,...,16} \draw (14,\x) node{$\bullet$};
 \draw[line width=1pt]  (14,0)--(14,16);
\foreach \x in {0,1,...,17} \draw (15,\x) node{$\bullet$};
 \draw[line width=1pt]  (15,0)--(15,17);
\foreach \x in {0,1,...,19} \draw (16,\x) node{$\bullet$};
 \draw[line width=1pt]  (16,0)--(16,19);
\foreach \x in {0,1,...,19} \draw (17,\x) node{$\bullet$};
 \draw[line width=1pt]  (17,0)--(17,19);
\foreach \x in {0,1,...,21} \draw (18,\x) node{$\bullet$};
 \draw[line width=1pt]  (18,0)--(18,21);
\foreach \x in {0,1,...,22} \draw (19,\x) node{$\bullet$};
 \draw[line width=1pt]  (19,0)--(19,22);
\foreach \x in {0,1,...,24} \draw (20,\x) node{$\bullet$};
 \draw[line width=1pt]  (20,0)--(20,24);
 \foreach \x in {0,1,...,24} \draw (21,\x) node{$\bullet$};
 \draw[line width=1pt]  (21,0)--(21,24);
\draw[color=red] (0,0) node{$\bullet$};  
\foreach \x in {0,1,...,5} \draw[color=red] (-1+\x,\x) node{$\bullet$};
 \draw[color=red, line width=1pt]  (-1,0)--(4,5);
\foreach \x in {0,1,...,12} \draw[color=red] (-2+\x,\x) node{$\bullet$};
 \draw[color=red, line width=1pt]  (-2,0)--(10,12);
\foreach \x in {0,1,...,17} \draw[color=red] (-3+\x,\x) node{$\bullet$};
 \draw[color=red, line width=1pt]  (-3,0)--(14,17);
\foreach \x in {0,1,...,24} \draw[color=red] (-4+\x,\x) node{$\bullet$};
 \draw[color=red, line width=1pt]  (-4,0)--(20,24);
\foreach \x in {0,1,...,26} \draw[color=red] (-5+\x,\x) node{$\bullet$};
 \draw[color=red, line width=1pt]  (-5,0)--(21.5,26.5);
 \foreach \x in {0,1,...,26} \draw[color=red] (-6+\x,\x) node{$\bullet$};
 \draw[color=red, line width=1pt]  (-6,0)--(20.5,26.5);
  \foreach \x in {1,2,...,26} \draw[color=red] (-7+\x,\x) node{$\bullet$};
 \draw[color=red, line width=1pt]  (-6.5,0.5)--(19.5,26.5);
   \foreach \x in {2,3,...,26} \draw[color=red] (-8+\x,\x) node{$\bullet$};
 \draw[color=red, line width=1pt]  (-6.5,1.5)--(18.5,26.5);
   \foreach \x in {3,4,...,26} \draw[color=red] (-9+\x,\x) node{$\bullet$};
 \draw[color=red, line width=1pt]  (-6.5,2.5)--(17.5,26.5);
    \foreach \x in {4,5,...,26} \draw[color=red] (-10+\x,\x) node{$\bullet$};
 \draw[color=red, line width=1pt]  (-6.5,3.5)--(16.5,26.5);
    \foreach \x in {5,6,...,26} \draw[color=red] (-11+\x,\x) node{$\bullet$};
 \draw[color=red, line width=1pt]  (-6.5,4.5)--(15.5,26.5);
     \foreach \x in {6,7,...,26} \draw[color=red] (-12+\x,\x) node{$\bullet$};
 \draw[color=red, line width=1pt]  (-6.5,5.5)--(14.5,26.5);
     \foreach \x in {7,8,...,26} \draw[color=red] (-13+\x,\x) node{$\bullet$};
 \draw[color=red, line width=1pt]  (-6.5,6.5)--(13.5,26.5);
     \foreach \x in {8,9,...,26} \draw[color=red] (-14+\x,\x) node{$\bullet$};
 \draw[color=red, line width=1pt]  (-6.5,7.5)--(12.5,26.5);
      \foreach \x in {9,10,...,26} \draw[color=red] (-15+\x,\x) node{$\bullet$};
 \draw[color=red, line width=1pt]  (-6.5,8.5)--(11.5,26.5);
       \foreach \x in {10,11,...,26} \draw[color=red] (-16+\x,\x) node{$\bullet$};
 \draw[color=red, line width=1pt]  (-6.5,9.5)--(10.5,26.5);
        \foreach \x in {11,12,...,26} \draw[color=red] (-17+\x,\x) node{$\bullet$};
 \draw[color=red, line width=1pt]  (-6.5,10.5)--(9.5,26.5);
         \foreach \x in {12,13,...,26} \draw[color=red] (-18+\x,\x) node{$\bullet$};
 \draw[color=red, line width=1pt]  (-6.5,11.5)--(8.5,26.5);
          \foreach \x in {13,14,...,26} \draw[color=red] (-19+\x,\x) node{$\bullet$};
 \draw[color=red, line width=1pt]  (-6.5,12.5)--(7.5,26.5);
           \foreach \x in {14,15,...,26} \draw[color=red] (-20+\x,\x) node{$\bullet$};
 \draw[color=red, line width=1pt]  (-6.5,13.5)--(6.5,26.5);
            \foreach \x in {15,16,...,26} \draw[color=red] (-21+\x,\x) node{$\bullet$};
 \draw[color=red, line width=1pt]  (-6.5,14.5)--(5.5,26.5);
\foreach \x in {16,17,...,26} \draw[color=red] (-22+\x,\x) node{$\bullet$};
 \draw[color=red, line width=1pt]  (-6.5,15.5)--(4.5,26.5);
\foreach \x in {17,18,...,26} \draw[color=red] (-23+\x,\x) node{$\bullet$};
 \draw[color=red, line width=1pt]  (-6.5,16.5)--(3.5,26.5);
 \foreach \x in {18,19,...,26} \draw[color=red] (-24+\x,\x) node{$\bullet$};
 \draw[color=red, line width=1pt]  (-6.5,17.5)--(2.5,26.5);
 \foreach \x in {19,20,...,26} \draw[color=red] (-25+\x,\x) node{$\bullet$};
 \draw[color=red, line width=1pt]  (-6.5,18.5)--(1.5,26.5);
 \foreach \x in {20,21,...,26} \draw[color=red] (-26+\x,\x) node{$\bullet$};
 \draw[color=red, line width=1pt]  (-6.5,19.5)--(0.5,26.5);
 \foreach \x in {21,22,...,26} \draw[color=red] (-27+\x,\x) node{$\bullet$};
 \draw[color=red, line width=1pt]  (-6.5,20.5)--(-0.5,26.5);
 \foreach \x in {22,23,...,26} \draw[color=red] (-28+\x,\x) node{$\bullet$};
 \draw[color=red, line width=1pt]  (-6.5,21.5)--(-1.5,26.5);
  \foreach \x in {23,24,25,26} \draw[color=red] (-29+\x,\x) node{$\bullet$};
 \draw[color=red, line width=1pt]  (-6.5,22.5)--(-2.5,26.5);
   \foreach \x in {24,25,26} \draw[color=red] (-30+\x,\x) node{$\bullet$};
 \draw[color=red, line width=1pt]  (-6.5,23.5)--(-3.5,26.5);
    \foreach \x in {25,26} \draw[color=red] (-31+\x,\x) node{$\bullet$};
 \draw[color=red, line width=1pt]  (-6.5,24.5)--(-4.5,26.5);
  \draw[color=red, line width=1pt]  (-6.5,25.5)--(-5.5,26.5);
 \draw[color=red] (-6,26) node{$\bullet$};
  \draw[color=blue, line width=1pt]  (1,1) circle[radius=0.3];
    \draw[color=blue, line width=1pt]  (5,6) circle[radius=0.3];
      \draw[color=blue, line width=1pt]  (7,8) circle[radius=0.3];
        \draw[color=blue, line width=1pt]  (11,13) circle[radius=0.3];
          \draw[color=blue, line width=1pt]  (15,18) circle[radius=0.3];
            \draw[color=blue, line width=1pt]  (17,20) circle[radius=0.3];
              \draw[color=blue, line width=1pt]  (21,25) circle[radius=0.3];
\draw (-5,0) node[below]{$-5$};
\draw (5,0) node[below]{$5$};
\draw (10,0) node[below]{$10$};
\draw (15,0) node[below]{$15$};
\draw (20,0) node[below]{$20$};
\draw (0,5) node[left]{$5$};
\draw (0,10) node[left]{$10$};
\draw (0,25) node[left]{$15$};
\draw (0,20) node[left]{$20$};
\draw (0,25) node[left]{$25$};
\draw (23,10) node{\large{$P_A$}};
\draw[color=red] (-8,13) node{\large{$P_B$}};
 \end{tikzpicture}
\end{equation}

The set $\{(\alpha,n) \mid \alpha \in \bZ \}$ in the above picture corresponds to
a $K$-basis of $F(n,n+1)$.
The set $\{(\alpha,n) \in P_A \mid \alpha \in \bZ \}$ corresponds to
a $K$-basis of $A(n,n+1)$.
The set $\{(\alpha,n) \in P_B \mid \alpha \in \bZ \}$ corresponds to a $K$-basis of $B(n,n+1)$.

\vspace{3mm}

\noindent
[\RMN{1}] \
First suppose ${\rm ch}(K)=2$.
We shall show that the condition (A4) in Theorem~\ref{chfree} is satisfied with $m=2$.
We have
\[
\xi^2 = (1-x)^4(1-x+vx)^{-2}
=(1-{x_{0,4}})(1-x_{0,2}+(x_{1,1})^2)^{-1}
=1+x_{0,2}-(x_{1,1})^2
\]
in $F_4$.
Here remark that $(x^iF)(x^jF)=x^{i+j}F$, $x_{\alpha,n} \in x^nF$ by (\ref{F}).
We know $F(2,4)=A(2,4)+B(2,4)$ by (\ref{pic}).
Take $f_A\in A(2,4)$ and $f_B\in B(2,4)$ such that $x_{0,2}-(x_{1,1})^2=f_A+f_B$.
Then we know $1+f_A \in A_4^\times$, $1+f_B \in \psi({B_4})^\times$ and
\[
(1+f_A)(1+f_B) = 1 + f_A+f_B = 1+x_{0,2}-(x_{1,1})^2
\]
in $F_4$.
Thus the condition (A4) in Theorem~\ref{chfree} is satisfied with $m=2$.
We know that ${\rm Cox}(Y)$ is Noetherian by Theorem~\ref{chfree}.

\vspace{3mm}

\noindent
[\RMN{2}] \
Next suppose ${\rm ch}(K)=3$.
We shall show that the condition (A4) in Theorem~\ref{chfree} is satisfied with $m=3$.
We have
\begin{align*}
\xi^3 = & (1-x)^6(1-x+vx)^{-3}
=(1-{x_{0,3}})^2(1-x_{0,3}+(x_{1,1})^3)^{-1} \\
= & (1-2x_{0,3})(1+x_{0,3}-(x_{1,1})^3)=1-x_{0,3}-(x_{1,1})^3
\end{align*}
in $F_6$.
We have $F(3,6)=A(3,6)+B(3,6)$ by (\ref{pic}).
Take $f'_A\in A(3,6)$ and $f'_B\in B(3,6)$ such that $-x_{0,3}-(x_{1,1})^3=f'_A+f'_B$.
Then we know $1+f'_A \in A_6^\times$, $1+f'_B \in \psi({B_6})^\times$ and
\[
(1+f'_A)(1+f'_B) = 1 + f'_A+f'_B = 1-x_{0,3}-(x_{1,1})^3
\]
in $F_6$.
Therefore the condition (A4) in Theorem~\ref{chfree} is satisfied with $m=3$.
We know that ${\rm Cox}(Y)$ is Noetherian by Theorem~\ref{chfree}.

\vspace{3mm}

\noindent
[\RMN{3}] \
Assume that the characteristic of $K$ is $p$, where $p$ is a prime number such that $p\ge 5$.
In the rest of this paper, we shall prove that ${\rm Cox}(Y)$ is not Noetherian.
It is enough to show that the condition (C3) in Theorem~\ref{chp} is not satisfied,
that is, we want to show 
\begin{equation}\label{wts}
\begin{array}{l}
\mbox{$H^0(\oo_Y(-\sigma jp^rC)|_{\sigma p^rC})= 0$
for any non-negative integer $r$} \\
\mbox{ and a positive integer $j$ such that $(j,p)=1$.
}
\end{array}
\end{equation}

Let $\chi(\ff)$ denotes the Euler characteristic of a coherent sheaf $\ff$ over $nC$,
that is,
\[
\chi(\ff) = \dim_KH^0(\ff) - \dim_KH^1(\ff) .
\]
By Claim~\ref{claim1} and (\ref{pic}), we know
\[
\oo_Y({-nC})/\oo_Y({-(n+1)C}) \simeq
\left\{
\begin{array}{ll}
\oo_{\bP^1_K} & (n \equiv 0 \mod 12) \\
\oo_{\bP^1_K}(-2) & (n \equiv 1,6,8 \mod 12) \\
\oo_{\bP^1_K}(-1) & (\mbox{otherwise}) .
\end{array}
\right.
\]
Therefore we have
\[
\chi(\oo_Y({-nC})/\oo_Y({-(n+1)C})) \simeq
\left\{
\begin{array}{cl}
1 & (n \equiv 0 \mod 12) \\
-1 & (n \equiv 1,6,8 \mod 12) \\
0 & (\mbox{otherwise}) .
\end{array}
\right.
\]
Since $\chi$ is an additive function for short exact sequences,
we know
\begin{equation}\label{chiofo}
\chi(\oo_Y(-\sigma jp^rC)|_{\sigma p^rC}) = -2p^r .
\end{equation}
%In order to prove (\ref{wts}), it is enough to show
%$\dim_KH^1(\oo_Y(-\sigma jp^rC)|_{\sigma p^rC})\le 2p^r$
%for any non-negative integer $r$
%and a positive integer $j$ such that $(j,p)=1$.

Put
\[
\maru{1}_d=x_{10d+1,12d+1}, \ \ \maru{2}_d=x_{10d+5,12d+6}, \ \ \maru{3}_d=x_{10d+7,12d+8}  
\]
and
\[
C_{p^r,j}=\{ \maru{i}_d
\mid {\rm i}=1,2,3;  \ d =  jp^r,   jp^r + 1, \ldots, (j+1) p^r-1 \} .
\]
By the exact sequence
\[
{\scriptstyle
0\longrightarrow A(\sigma jp^r, \sigma (j+1)p^r)+ B(\sigma jp^r, \sigma (j+1)p^r)
\longrightarrow F(\sigma jp^r, \sigma (j+1)p^r)
\longrightarrow H^1(\oo_Y(-\sigma jp^rC)|_{\sigma p^rC})
\longrightarrow 0 ,
}
\]
$C_{p^r,j}$ spans $H^1(\oo_Y(-\sigma jp^rC)|_{\sigma p^rC})$ as a $K$-vector space.

We define the subset $D_{p^r,j}$ of $C_{p^r,j}$ as
\[
D_{p^r,j} = \left\{ x_{\alpha,n} \in C_{p^r,j} \ \left| \
\mbox{$x_{\alpha,n}\equiv 0$ modulo $A(\sigma jp^r, n+1) + B(\sigma jp^r, n+1)$ in $F(\sigma jp^r, n+1)$} 
\right\} \right. .
\]

\begin{Lemma}\label{ablemma}
The set $C_{p^r,j}\setminus D_{p^r,j}$ is a $K$-basis of  $H^1(\oo_Y(-\sigma jp^rC)|_{\sigma p^rC})$.
\end{Lemma}

\proof
First we shall prove that $C_{p^r,j}\setminus D_{p^r,j}$ spans $H^1(\oo_Y(-\sigma jp^rC)|_{\sigma p^rC})$ as a $K$-vector space.
Let $\mathop{<}C_{p^r,j}\setminus D_{p^r,j}\mathop{>}_K$ be the $K$-vector subspace of $F(\sigma jp^r, \sigma (j+1)p^r)$ spanned by $C_{p^r,j}\setminus D_{p^r,j}$.
We shall prove 
\begin{equation}\label{nind}
A(\sigma jp^r, \sigma (j+1)p^r) + B(\sigma jp^r, \sigma (j+1)p^r) + \mathop{<}C_{p^r,j}\setminus D_{p^r,j}\mathop{>}_K
 + F(n, \sigma (j+1)p^r) = F(\sigma jp^r, \sigma (j+1)p^r)
\end{equation}
for $n = \sigma jp^r, \sigma jp^r+1, \ldots, \sigma (j+1)p^r$
by induction on $n$.
It is obvious in the case $n=\sigma jp^r$.
It is enough to show
\[
A(\sigma jp^r, \sigma (j+1)p^r) + B(\sigma jp^r, \sigma (j+1)p^r)  + \mathop{<}C_{p^r,j}\setminus D_{p^r,j}\mathop{>}_K
 + F(n+1, \sigma (j+1)p^r) \supset  F(n, \sigma (j+1)p^r)
\]
for $n = \sigma jp^r, \sigma jp^r+1, \ldots, \sigma (j+1)p^r-1$.
We have only to show that each $x_{\alpha,n}$ ($\alpha \in \bZ$) is contained in the left-hand side.
If $(\alpha,n)\in P_A$, then $x_{\alpha,n}$ is contained in $A(\sigma jp^r, \sigma (j+1)p^r)$.
If $(\alpha,n)\in P_B$, then $x_{\alpha,n}$ is contained in $B(\sigma jp^r, \sigma (j+1)p^r)+F(n+1, \sigma (j+1)p^r)$ by (4.15) in \cite{K43}.
If $x_{\alpha,n}\in C_{p^r,j}\setminus D_{p^r,j}$, then $x_{\alpha,n}$ is contained in $\mathop{<}C_{p^r,j}\setminus D_{p^r,j}\mathop{>}_K$.
If $x_{\alpha,n}\in D_{p^r,j}$, then $x_{\alpha,n}$ is contained in the left-hand side by definition of $D_{p^r,j}$.
Therefore (\ref{nind}) is satisfied for $n = \sigma (j+1)p^r$.
Hence $C_{p^r,j}\setminus D_{p^r,j}$ spans $H^1(\oo_Y(-\sigma jp^rC)|_{\sigma p^rC})$ as a $K$-vector space.

Next we shall prove that  $C_{p^r,j}\setminus D_{p^r,j}$ are linearly independent in  $H^1(\oo_Y(-\sigma jp^rC)|_{\sigma p^rC})$.
Assume the contrary.
There exist $x_{\alpha_1,n_1}$, $x_{\alpha_2,n_2}$, \ldots, $x_{\alpha_k,n_k} \in C_{p^r,j}\setminus D_{p^r,j}$ and $c_1, \ldots, c_k \in K\setminus \{ 0 \}$ such that
\[
c_1x_{\alpha_1,n_1} + c_2x_{\alpha_2,n_2}+ \cdots+c_kx_{\alpha_k,n_k} =0
\]
in $H^1(\oo_Y(-\sigma jp^rC)|_{\sigma p^rC})$.
Suppose $n_1<n_2<\cdots<n_k$.
Then we obtain $x_{\alpha_1,n_1} \in D_{p^r,j}$.
It is a contradiction.
\qed

\begin{Remark}\label{abrem}
\begin{rm}
Remark ${}^\#C_{p^r,j} = 3p^r$.

Here assume ${}^\#D_{p^r,j} \ge p^r$.
Then we know 
\[
\dim_KH^1(\oo_Y(-\sigma jp^rC)|_{\sigma p^rC})\le 2p^r 
\]
by Lemma~\ref{ablemma}.
On the other hand we know its Euler characteristic (\ref{chiofo}).
Therefore we have ${}^\#D_{p^r,j} = p^r$, $\dim_KH^1(\oo_Y(-\sigma jp^rC)|_{\sigma p^rC})= 2p^r$ and $H^0(\oo_Y(-\sigma jp^rC)|_{\sigma p^rC})= 0$.

Therefore, in order to prove (\ref{wts}), it is enough to show
\begin{equation}\label{ets}
\begin{array}{l}
\mbox{${}^\#D_{p^r,j} \ge p^r$
for any non-negative integer $r$} \\
\mbox{ and a positive integer $j$ such that $(j,p)=1$.
}
\end{array}
\end{equation}
\end{rm}
\end{Remark}

In the rest of this paper, we shall prove (\ref{ets})
for each prime number $p$ such that $p\ge 5$.

By (4.5) and (4.6) in \cite{K43}, $x_{\alpha,n}$ and $z_{\alpha,n}$ are contained in $x^nF$.
By definition (\ref{xalphan}),
we know
\[
x_{\alpha,n}x_{\alpha',n'}=x_{\alpha+\alpha',n+n'}
\]
if either $\alpha$ or $\alpha'$ is even.
In particular, we have
\[
x_{\alpha,n}x^m = x_{\alpha,n}(x_{0,1})^m = x_{\alpha,n+m} .
\]
By (\ref{zalphan}),
we know
\[
z_{\alpha,n}z_{\alpha',n'}=z_{\alpha+\alpha',n+n'}
\]
if either $\alpha-n$ or $\alpha'-n'$ is even.
We have
\[
x_{\alpha,n}=z_{\alpha,n}
\]
in $F(n,n+1)$ by  (4.15) in \cite{K43}.
We have
\begin{equation}\label{xnoseki}
x_{\alpha_1,n_1}x_{\alpha_2,n_2} \cdots x_{\alpha_t,n_t}=x_{\alpha_1+\cdots+\alpha_t,n_1+\cdots+n_t}
=z_{\alpha_1+\cdots+\alpha_t,n_1+\cdots+n_t}
=z_{\alpha_1,n_1}z_{\alpha_2,n_2} \cdots z_{\alpha_t,n_t}
\end{equation}
in $F(n_1+\cdots+n_t,n_1+\cdots+n_t+1)$ by  (4.8) in \cite{K43}.

We have
\begin{equation}\label{ztox}
\mbox{$z_{10k,12k} = x_{10k,12k} \xi^{-6k}$ \ where \ $\xi=(1-x)^2(1-x+vx)^{-1}$}
\end{equation}
by definition (\ref{zalphan}).

Suppose that $k$ is a positive integer such that $(k,p)=1$.
We have
\begin{align*}
z_{10k,12k}-x_{10k,12k} & = x_{10k,12k}(\xi^{-6k}-1) \\
& = x_{10k,12k}((1-x)^{-12k}(1-x+vx)^{6k}-1) \\
& \equiv x_{10k,12k}(6kx+6kvx) \\
& \equiv (6k)z_{10k,12k+1} + (6k)x_{10k+1,12k+1} 
\end{align*}
modulo $x^{12k+2}F$.

Here suppose  $jp^r \le k < (j+1)p^r$.
Remark that $x_{10k,12k} \in A(12jp^r,12k+2)$,
and $z_{10k,12k}, z_{10k,12k+1}\in B(12jp^r,12k+2)$.
Recall that the characteristic of the base field $K$ is not $2$ or $3$.
Since $6k \neq 0$, we know that 
\begin{equation}\label{tagainiso1} 
x_{10k,12k}vx = x_{10k+1,12k+1}\equiv 0
\end{equation}
modulo $A(12jp^r,12k+2)+B(12jp^r,12k+2)$
in $F(12jp^r,12k+2)$.
Thus we know
\begin{equation}\label{tagainiso2}
\mbox{$\maru{1}_k$ is in $D_{p^r,j}$ if $jp^r \le k < (j+1)p^r$ and $(k,p)=1$.}
\end{equation}

The following lemma will be frequently used later.

\begin{Lemma}\label{omeganoseki}
Let $k$ be a positive integer.
\begin{enumerate}
\item
Assume that $\alpha$ is even.
Then
\begin{align*}
x_{\alpha,n}w^k = & \sum_{q=0}^{k-1}\left\{
{k+q-1\choose 2q}(1-x)^{k-q}x_{\alpha+2q,n+2q} +
{k+q\choose 2q+1}(1-x)^{k-q-1}x_{\alpha+2q+1,n+2q+1} 
\right\} \\
= & (1-x)^kx_{\alpha,n} + {k \choose 1}(1-x)^{k-1}x_{\alpha+1,n+1} \\
&+ {k \choose 2}(1-x)^{k-1}x_{\alpha+2,n+2}
 + {k+1 \choose 3}(1-x)^{k-2}x_{\alpha+3,n+3}  \\
 &+ {k+1 \choose 4}(1-x)^{k-2}x_{\alpha+4,n+4} 
 + {k+2 \choose 5}(1-x)^{k-3}x_{\alpha+5,n+5}  \\ &+ {k+2 \choose 6}(1-x)^{k-3}x_{\alpha+6,n+6} + {k+3 \choose 7}(1-x)^{k-4}x_{\alpha+7,n+7}+\cdots
\end{align*}
\item
Assume that $\alpha$ is odd.
Then
\begin{align*}
x_{\alpha,n}w^k = & (1-x)^kx_{\alpha,n} \\
& + \sum_{q=0}^{k-1}\left\{
{k+q\choose 2q+1}(1-x)^{k-q}x_{\alpha+2q+1,n+2q+1} +
{k+q+1\choose 2q+2}(1-x)^{k-q-1}x_{\alpha+2q+2,n+2q+2} 
\right\} \\
= & (1-x)^kx_{\alpha,n} \\
&  + {k \choose 1}(1-x)^{k}x_{\alpha+1,n+1}
+ {k+1 \choose 2}(1-x)^{k-1}x_{\alpha+2,n+2} \\
& + {k+1 \choose 3}(1-x)^{k-1}x_{\alpha+3,n+3} + {k+2 \choose 4}(1-x)^{k-2}x_{\alpha+4,n+4}  \\
& + {k+2 \choose 5}(1-x)^{k-2}x_{\alpha+5,n+5}  + {k+3 \choose 6}(1-x)^{k-3}x_{\alpha+6,n+6} + \cdots
\end{align*}
\end{enumerate}
\end{Lemma}

\proof
First assume that $\alpha$ is even.
By definition (\ref{xalphan}), we have
\[
x_{\alpha+2q,n+2q}=x_{\alpha,n}(vx)^{2q}w^{-q}, \ \ 
x_{\alpha+2q+1,n+2q+1}=x_{\alpha,n}(vx)^{2q+1}w^{-q} .
\]
Put
\begin{align*}
A_q & = \sum_{i=2q}^{k+q-1}{k+q-1 \choose i}(1-x)^{k-i+q-1}(vx)^iw^{-(q-1)} , \\
B_q & = {k+q-1 \choose 2q}(1-x)^{k-q}(vx)^{2q}w^{-q}
+ {k+q \choose 2q+1}(1-x)^{k-q-1}(vx)^{2q+1}w^{-q} .
\end{align*}
It is enough to show
\begin{equation}\label{w^k}
w^k = \sum_{q=0}^{k-1}B_q .
\end{equation}
We have
\begin{align}\label{4.8}
 A_q 
= & \sum_{i=2q}^{k+q-1}{k+q-1 \choose i}(1-x)^{k-i+q-1}(vx)^iw^{-q}(1-x+vx) \nonumber \\
= & \sum_{i=2q}^{k+q-1}{k+q-1 \choose i}(1-x)^{k-i+q}(vx)^iw^{-q}
+ \sum_{j=2q+1}^{k+q}{k+q-1 \choose j-1}(1-x)^{k-j+q}(vx)^jw^{-q}  \nonumber\\
= & {k+q-1 \choose 2q}(1-x)^{k-q}(vx)^{2q}w^{-q}
+ \sum_{i=2q+1}^{k+q}{k+q \choose i}(1-x)^{k-i+q}(vx)^iw^{-q}  \nonumber\\
= & {k+q-1 \choose 2q}(1-x)^{k-q}(vx)^{2q}w^{-q}
+ {k+q \choose 2q+1}(1-x)^{k-q-1}(vx)^{2q+1}w^{-q}  \nonumber\\
& + \sum_{i=2(q+1)}^{k+q}{k+q \choose i}(1-x)^{k-i+q}(vx)^iw^{-q} . \nonumber\\
= & B_q + A_{q+1} \nonumber
\end{align}
Using this formula several times, we have
\begin{align*}
w^k = & (1-x+vx)^k = \sum_{i=0}^k{k \choose i}(1-x)^{k-i}(vx)^i = B_0+A_1 = B_0+B_1+A_2 = \cdots \\
= & \sum_{q=0}^{k-2}B_q+A_{k-1} = \sum_{q=0}^{k-1}B_q
\end{align*}
since $A_{k-1}=B_{k-1}$.
We have proved (1).

Next assume that $\alpha$ is odd.
By definition (\ref{xalphan}), we have
\[
x_{\alpha+2q,n+2q}=x_{\alpha,n}(vx)^{2q}w^{-q}, \ \ 
x_{\alpha+2q+1,n+2q+1}=x_{\alpha,n}(vx)^{2q+1}w^{-q-1} .
\]
It is enough to show
\begin{align*}
w^k = & (1-x)^k \\
& + \sum_{q=0}^{k-1}\left\{
{k+q\choose 2q+1}(1-x)^{k-q}(vx)^{2q+1}w^{-q-1} +
{k+q+1\choose 2q+2}(1-x)^{k-q-1}(vx)^{2q+2}w^{-q-1}
\right\}
\end{align*}
By (\ref{w^k}), we obtain
\begin{align*}
w^k = & 
\sum_{q=0}^{k-1}{k+q-1 \choose 2q}(1-x)^{k-q}(vx)^{2q}w^{-q}
+\sum_{q=0}^{k-1}{k+q \choose 2q+1}(1-x)^{k-q-1}(vx)^{2q+1}w^{-q}
\\
%= & (1-x)^k+\sum_{q=1}^{k-1}{k+q-1 \choose 2q}(1-x)^{k-q}(vx)^{2q}w^{-q} \\
%& +\sum_{q=0}^{k-1}{k+q \choose 2q+1}(1-x)^{k-q}(vx)^{2q+1}w^{-q-1}
%+ \sum_{q=0}^{k-1}{k+q \choose 2q+1}(1-x)^{k-q-1}(vx)^{2q+2}w^{-q-1} \\
= & (1-x)^k+\sum_{q'=0}^{k-2}{k+q' \choose 2q'+2}(1-x)^{k-q'-1}(vx)^{2q'+2}w^{-q'-1} \\
& +\sum_{q=0}^{k-1}{k+q \choose 2q+1}(1-x)^{k-q}(vx)^{2q+1}w^{-q-1}
+ \sum_{q=0}^{k-1}{k+q \choose 2q+1}(1-x)^{k-q-1}(vx)^{2q+2}w^{-q-1} \\
= & (1-x)^k\\
& + \sum_{q=0}^{k-1}{k+q \choose 2q+1}(1-x)^{k-q}(vx)^{2q+1}w^{-q-1}
+ \sum_{q=0}^{k-1}{k+q+1 \choose 2q+2}(1-x)^{k-q-1}(vx)^{2q+2}w^{-q-1} .
\end{align*}
\qed

\vspace{3mm}

Remark that, by this lemma, we can rewite $x_{\alpha,n}w^k$ into a $K$-linear combination of $x_{\beta,m}$'s since $x_{\beta,m}x^\ell = x_{\beta,m+\ell}$.

\vspace{3mm}

In the rest, we shall prove Example~\ref{rei} (2) by dividing into some cases.

\vspace{3mm}

\noindent
[\RMN{3}-1] \
Assume that $p$ is $5$.
Taking $p$th power of the equation (\ref{tagainiso1}), we obtain
\[
(x_{10k,12k}vx)^p \equiv 0
\]
modulo $A(12kp,12kp+2p)+B(12kp,12kp+2p)$ in $F(12kp,12kp+2p)$.
We have
\begin{align*}
& (x_{10k,12k} vx)^p = x_{10kp,12kp} (vx)^5=x_{10kp,12kp}v^5w^{-2}x^5w^2 \\
= & x_{10kp+5,12kp+5}(1-x+vx)^2 = x_{10kp+5,12kp+5}(1-2x+2vx+x^2 -2vx^2 +v^2x^2) \\
 \equiv & x_{10kp+5,12kp+5} -2 x_{10kp+5,12kp+6} + 2 x_{10kp+6,12kp+6}
\end{align*}
modulo $x^{12kp+7}F$.
Since $x_{10kp+5,12kp+5}, x_{10kp+6,12kp+6} \in A(12kp,12kp+7)$,
we have 
\[
x_{10kp+5,12kp+6}\equiv 0
\]
modulo $A(12kp,12kp+7)+B(12kp,12kp+7)$ in $F(12kp,12kp+7)$.
Taking $p^{h-1}$th power of it for $h>0$, we have
\[
x_{10kp^h+5p^{h-1},12kp^h+6p^{h-1}}\equiv 0
\]
modulo $A(12kp^h,12kp^h+6p^{h-1}+1)+B(12kp^h,12kp^h+6p^{h-1}+1)$ in $F(12kp^h,12kp^h+6p^{h-1}+1)$ by (\ref{xnoseki}).
%Here we have $x_{10p^hk,12p^hk}, (x_{10pk+5,12pk+5})^{p^{h-1}}, x_{10p^hk+6,12p^hk+6}  \in A(12p^hk, 12p^hk+6p^{h-1}+1)$ and $z_{10p^hk,12p^hk}, (z_{10k,12k+1})^{p^h} \in B(12p^hk, 12p^hk+6p^{h-1}+1)$.
Here remark 
\begin{align*}
10kp^h+5p^{h-1} & = 10\left( kp^h + \frac{p-1}{2}(p^{h-2} + p^{h-3}+ \cdots + p+1)  \right) +5 , \\
12kp^h+6p^{h-1} & = 12\left( kp^h + \frac{p-1}{2}(p^{h-2} + p^{h-3}+ \cdots + p+1)  \right) +6 .
\end{align*}
Therefore we know
\begin{equation}\label{pnobaisuu}
\mbox{$\maru{2}_{kp^h + 2(p^{h-2} + p^{h-3}+ \cdots + p+1) }$ is in $D_{p^r,j}$ if $r \ge h>0$, $jp^r \le kp^h < (j+1)p^r$ and $(k,p)=1$.}
\end{equation}
Then $D_{p^r,j}$ contains $p^r-p^{r-1}$ elements of the form $\maru{1}_d$ as in (\ref{tagainiso2}) and $p^{r-1}$ elements of the form $\maru{2}_d$ as in (\ref{pnobaisuu}).
Therefore $D_{p^r,j}$ contains at least $p^r$ elements.
We know that ${\rm Cox}(Y)$ is not Noetherian by (\ref{ets}).

\vspace{3mm}

\noindent
[\RMN{3}-2] \
Assume that $p=10f+3$, where $f >0$.
let $k$ be a positive integer such that $(k,p)=1$.
Put $e = kp+f$.
First we shall prove
\begin{equation}\label{aim10f+3}
x_{10e+5,12e+6}\equiv 0
\end{equation}
modulo $A(12kp,12e+7)+B(12kp,12e+7)$
in $F(12kp,12e+7)$.

Taking the $p$th power of  (\ref{tagainiso1}), we obtain
\begin{equation}\label{10f3-1}
x_{10kp,12kp}(vx)^p \equiv 0 
\end{equation}
modulo $A(12kp,12e+7)+B(12kp,12e+7)$
in $F(12kp,12e+7)$
since $12kp+2p > 12(kp+f)+7=12e+7$.
We have
\begin{align}
& x_{10kp,12kp}(vx)^p \label{10f3-2}
=  x_{10kp,12kp} x_{p,p}w^{5f+1}  \\
%= &  x_{10kp+p,12kp+p}
%\left(
%(1-x)^{5f+1} + (5f+1)(1-x)^{5f+1}vxw^{-1} + {5f+2 \choose 2}(1-x)^{5f}(vx)^2w^{-1} 
%+\cdots 
%\right)   \nonumber 
%\\ 
\nonumber 
= & (1-x)^{5f+1}x_{10kp+p,12kp+p} + (5f+1)(1-x)^{5f+1}x_{10kp+p+1,12kp+p+1} \\
& +
{5f+2 \choose 2}(1-x)^{5f}x_{10kp+p+2,12kp+p+2}+\cdots 
\nonumber 
\end{align}
by Lemma~\ref{omeganoseki} (2).
Here recall that $x_{\alpha,n}x^q= x_{\alpha,n+q}$.
If $x_{\alpha,n}$ appears in (\ref{10f3-2}),
$(\alpha,n)$ safisfies
\begin{equation}\label{10f3-3}
\left\{
\begin{array}{l}
\alpha \ge 10kp+p \\
n \ge \alpha + 2kp
\end{array}
\right.
\end{equation}
Remark that, if $n \ge 10e+7$, then $x_{\alpha,n}=0$ in $F(12kp,12e+7).$
\begin{equation}\label{pic2}
\begin{tikzpicture}[xscale = 0.45, yscale = 0.45] 
\draw[color=gray] (0.5,-0.5) grid(14.5,13.5);
\draw[line width=0.8pt] (0.5,-0.5) grid[step=5](14.5,13.5);
%\draw[->,>=stealth,semithick] (-6,0)--(22,0)node[below]{}; %x軸
%\draw[->,>=stealth,thick] (0,-2)--(0,15)node[right]{}; %y軸
\draw[color=green, line width=1pt]  (0.5,12)--(14.5,12);
\draw (0.5,12)node[left]{$12e+7$}; 
\draw (0.5,8)node[left]{$12e+3$}; 
\draw (0.5,6)node[left]{$12kp+p$}; 
\draw[color=blue, line width=1pt]  (10,11) circle[radius=0.3];
\foreach \x in {0,1,...,5} \draw (8+\x,6+\x) node{$\bullet$};
\foreach \x in {0,1,...,4} \draw (8+\x,7+\x) node{$\bullet$};
\foreach \x in {0,1,...,3} \draw (8+\x,8+\x) node{$\bullet$};
\foreach \x in {0,1,2} \draw[color=red] (8,9+\x) node{$\bullet$};
\foreach \x in {0,1} \draw[color=red] (9,10+\x) node{$\bullet$};
\draw[line width=1pt] (8,-0.5)--(8,0);
\draw (8,-0.5)node[below]{$10kp+p$}; 
 \end{tikzpicture}
\end{equation}

Remark $10kp+p=10e+3$.
If $(\alpha,n)$ in the area (\ref{10f3-3}) satisfies $n \le \alpha+2e$,
then $x_{\alpha,n}\in A(12kp,12e+7)$.
We know that $x_{10e+3,12e+5}$, $x_{10e+3,12e+6}$, 
$x_{10e+4,12e+6}$ are equivalent to 0 modulo $B(12kp,12e+7)$
in $F(12kp,12e+7)$ by (4.15) in \cite{K43}.
Since $(12e+4)-(12kp+p)=2f+1$, we have
\begin{align}
& x_{10kp,12kp}(vx)^p \label{10f3-4} \\
\equiv & -{5f+1\choose 2f+1} x_{10e+3,12e+4} -(5f+1){5f+1\choose 2f+1} x_{10e+4,12e+5}
-{5f+2\choose 2}{5f\choose 2f+1}x_{10e+5,12e+6} \nonumber \\
=&  -{5f+1\choose 2f+1}\left(
x_{10e+3,12e+4} +(5f+1) x_{10e+4,12e+5}
+\frac{(5f+2)3f}{2}x_{10e+5,12e+6} 
\right) \nonumber 
\end{align}
modulo $A(12kp,12e+7)+B(12kp,12e+7)$
in $F(12kp,12e+7)$.
Here $-{5f+1\choose 2f+1}\neq 0$ in $K$.
By (\ref{10f3-1}), (\ref{10f3-4}),
we know 
\begin{equation}\label{10f3-5}
x_{10e+3,12e+4} +(5f+1) x_{10e+4,12e+5}
+\frac{(5f+2)3f}{2}x_{10e+5,12e+6} \equiv 0
\end{equation}
modulo $A(12kp,12e+7)+B(12kp,12e+7)$
in $F(12kp,12e+7)$.
Furthermore we have 
\begin{align}
\label{10f3-6}
0 \equiv & z_{10e+3,12e+4} = v^{10e+3}w^{e+1}x^{12e+4}(1+x+x^2+\cdots)^{12e+4} \\
\equiv & x_{10e+3,12e+4}w^{6e+2} \nonumber \\
\equiv & 
x_{10e+3,12e+4} +(6e+2) x_{10e+4,12e+5}
+{6e+3 \choose 2}x_{10e+5,12e+6} \nonumber
\end{align}
modulo $A(12kp,12e+7)+B(12kp,12e+7)$
in $F(12kp,12e+7)$
by Lemma~\ref{omeganoseki} (2).
We have 
\begin{align}
\label{10f3-7}
0 \equiv & z_{10e+4,12e+5} = v^{10e+4}w^{e+1}x^{12e+5}(1+x+x^2+\cdots)^{12e+5} \\
\equiv & x_{10e+4,12e+5}w^{6e+3} \nonumber \\
\equiv & 
x_{10e+4,12e+5}
+(6e+3)x_{10e+5,12e+6} \nonumber
\end{align}
modulo $A(12kp,12e+7)+B(12kp,12e+7)$
in $F(12kp,12e+7)$
by Lemma~\ref{omeganoseki} (1).
Here remark $e=kp+f \equiv f \mod p$.
By (\ref{10f3-5}),  (\ref{10f3-6}),  (\ref{10f3-7}),
we obtain 
\[
\frac{3f(3f+2)}{2} x_{10e+5,12e+6}\equiv 0
\]
modulo $A(12kp,12e+7)+B(12kp,12e+7)$
in $F(12kp,12e+7)$.
Since $\frac{3f(3f+2)}{2} \neq 0$ in $K$, we obtain (\ref{aim10f+3}).

For $h>0$, taking the $p^{h-1}$th power of the above equation,
we obtain
\[
x_{10ep^{h-1}+5p^{h-1},12ep^{h-1}+6p^{h-1}}\equiv 0
\]
modulo $A(12kp^h,12ep^{h-1}+6p^{h-1}+1)+B(12kp^h,12ep^{h-1}+6p^{h-1}+1)$
in $F(12kp^h,12ep^{h-1}+6p^{h-1}+1)$ by (\ref{xnoseki}).
Here remark 
\begin{align*}
10ep^{h-1}+5p^{h-1} = & 10kp^h+10fp^{h-1}+5p^{h-1} \\
= & 10\left( kp^h + fp^{h-1}+\frac{p-1}{2}(p^{h-2} + p^{h-3}+ \cdots + p+1)  \right) +5 , \\
12ep^{h-1}+6p^{h-1} = & 12kp^h+12fp^{h-1}+6p^{h-1} \\
= & 12\left( kp^h + fp^{h-1}+ \frac{p-1}{2}(p^{h-2} + p^{h-3}+ \cdots + p+1)  \right) +6 .
\end{align*}
Therefore we know
\begin{equation}\label{pnobaisuu3}
\mbox{$\maru{2}_{kp^h + fp^{h-1}+(5f+1)(p^{h-2} + p^{h-3}+ \cdots + p+1) }$ is in $D_{p^r,j}$ if $r \ge h>0$, $jp^r \le kp^h < (j+1)p^r$ and $(k,p)=1$.}
\end{equation}
Then $D_{p^r,j}$ contains $p^r-p^{r-1}$ elements of the form $\maru{1}_d$ as in (\ref{tagainiso2}) and $p^{r-1}$ elements of the form $\maru{2}_d$ as in (\ref{pnobaisuu3}).
Therefore $D_{p^r,j}$ contains at least $p^r$ elements.
We know that ${\rm Cox}(Y)$ is not Noetherian by (\ref{ets}).

\vspace{3mm}

\noindent
[\RMN{3}-3] \
Assume that $p=10f+1$, where $f >0$.

Let $k$ be a positive integer such that $(k,p)=1$.
Put $e=kp+f$.

First we shall prove that
\begin{equation}\label{10f1-0}
x_{10e+5,12e+6} \equiv 0
\end{equation}
modulo $A(12kp,12e+7)+B(12kp,12e+7)$ in $F(12kp,12e+7)$.

Consider the triangle $T$ with three vertices $(10e,12kp+10f)$,
$(10e,12e+6)$, $(10e+6+2f,12e+6)$.
\begin{equation}\label{pic5}
\begin{tikzpicture}[xscale = 0.45, yscale = 0.45] 
\draw[color=gray] (0.5,-0.5) grid(14.5,13.5);
\draw[line width=0.8pt] (0.5,-0.5) grid[step=5](14.5,13.5);
%\draw[->,>=stealth,semithick] (-6,0)--(22,0)node[below]{}; %x軸
%\draw[->,>=stealth,thick] (0,-2)--(0,15)node[right]{}; %y軸
\fill (5,5) circle[radius=0.3];
\draw[color=green, line width=1pt]  (0.5,12)--(14.5,12);
\draw (0.5,12)node[left]{$12e+7$}; 
\draw (0.5,11)node[left]{$12e+6$}; 
\draw (0.5,5)node[left]{$12e$}; 
\draw (0.5,3)node[left]{$12kp+10f$}; 
\draw[color=blue, line width=1pt]  (10,11) circle[radius=0.3];
\draw[color=blue, line width=1pt]  (6,6) circle[radius=0.3];
\draw (5,5) node{$\bullet$};
\foreach \x in {0,1,...,8} \draw (5+\x,3+\x) node{$\bullet$};
\foreach \x in {0,1,...,7} \draw (5+\x,4+\x) node{$\bullet$};
\foreach \x in {0,1,...,4} \draw (7+\x,7+\x) node{$\bullet$};
\draw[color=red] (5,11) node{$\bullet$};
\draw[color=red] (5,5) node{$\bullet$};
\foreach \x in {0,1} \draw[color=red] (5+\x,10+\x) node{$\bullet$};
\foreach \x in {0,1,2} \draw[color=red] (5+\x,9+\x) node{$\bullet$};
\foreach \x in {0,1,...,3} \draw[color=red] (5+\x,8+\x) node{$\bullet$};
\foreach \x in {0,1,...,4} \draw[color=red] (5+\x,7+\x) node{$\bullet$};
\foreach \x in {0,1,...,4} \draw[color=red] (5+\x,6+\x) node{$\bullet$};
\draw[line width=1pt] (13,-0.5)--(13,0);
\draw (5,-0.5)node[below]{$10e$}; 
\draw (13,-0.5)node[below]{$10e+6+2f$}; 
 \end{tikzpicture}
\end{equation}
For $c_0,c_1,\ldots,c_6 \in K$, we put
\begin{align*}
 [c_0,c_1,\ldots,c_6]
= & c_0x_{10e+1,12e+1}+c_1x_{10e,12e+1}+c_2x_{10e+1,12e+2}+c_3x_{10e+2,12e+3} \\
& +c_4x_{10e+3,12e+4}+c_5x_{10e+4,12e+5}+c_6x_{10e+5,12e+6} .
\end{align*}
If $(\alpha,n)$ is in the triangle $T$ such that $n \le \alpha + 2e$,
then $x_{\alpha,n}$ is in $A(12kp+10f,12e+7)$ except for $x_{10e+1,12e+1}$.
If $(\alpha,n)$ is in the triangle $T$ such that $n \ge \alpha + 2e+2$,
then $x_{\alpha,n}$ is in $B(12kp+10f,12e+7)$.
Therefore any $K$-linear combination of $x_{\alpha,n}$'s in the triangle $T$ is equivalent to some  $[c_0,c_1,\ldots,c_6]$ 
modulo $A(12kp+10f,12e+7)+B(12kp+10f,12e+7)$ in $F(12kp+10f,12e+7)$.

Taking the $p$th power of (\ref{tagainiso1}), we obtain 
\[
x_{10kp,12kp}(vx)^p \equiv 0
\]
modulo $A(12kp,12kp+7)+B(12kp,12kp+7)$ in $F(12kp,12kp+7)$
since $12kp+2p=12kp+20f+2>12kp+12f+7=12e+7$.
Here remark that $x_{10kp,12kp}(vx)^p \mod x^{12e+7}F$ is a $K$-linear combination of
$x_{\alpha,n}$'s in the triangle $T$
since $(10kp+p,12kp+p)=(10e+1,12kp+10f+1) \in T$.
We have
\begin{align}
\label{10f1-2}
& x_{10kp,12kp}(vx)^p=x_{10kp+p,12kp+p}w^{5f} = x_{10e+1,12e+1-2f}w^{5f} \\
= & 
 (1-x)^{5f}x_{10e+1,12e+1-2f} + {5f \choose 1}(1-x)^{5f}x_{10e+2,12e+2-2f}+ {5f+1 \choose 2}(1-x)^{5f-1}x_{10e+3,12e+3-2f}  \nonumber
\\
& + {5f+1 \choose 3}(1-x)^{5f-1}x_{10e+4,12e+4-2f} + 
 {5f+2 \choose 4}(1-x)^{5f-2}x_{10e+5,12e+5-2f} + 
\cdots  \nonumber \\
\equiv & \left[
{\textstyle
{5f \choose 2f}, 0, -{5f \choose 2f+1}, -5f{5f \choose 2f+1}, 
-{5f+1 \choose 2}{5f-1 \choose 2f+1}, -{5f+1 \choose 3}{5f-1 \choose 2f+1}, 
-{5f+2 \choose 4}{5f-2 \choose 2f+1}
}
\right] \nonumber
\end{align}
modulo $A(12kp+10f,12e+7)+B(12kp+10f,12e+7)$ in $F(12kp+10f,12e+7)$
by Lemma~\ref{omeganoseki} (2).

By (\ref{ztox}) we obtain
\begin{align}
\label{10f1-3}
0 \equiv & z_{10e,12e}=x_{10e,12e}\xi^{-6e} =
(1-x)^{-12e} x_{10e,12e} w^{6e}\\
= & 
(1-x)^{-6e}x_{10e,12e} + {6e \choose 1}(1-x)^{-6e-1}x_{10e+1,12e+1}+ {6e \choose 2}(1-x)^{-6e-1}x_{10e+2,12e+2}  \nonumber
\\
& + {6e+1 \choose 3}(1-x)^{-6e-2}x_{10e+3,12e+3} + 
 {6e+1 \choose 4}(1-x)^{-6e-2}x_{10e+4,12e+4} \nonumber \\
 &+  {6e+2 \choose 5}(1-x)^{-6e-3}x_{10e+5,12e+5} +
\cdots  \nonumber \\
\equiv & \left[
{\textstyle
6e, 6e, 6e(6e+1), {6e \choose 2}(6e+1), 
{6e+1 \choose 3}(6e+2), {6e+1 \choose 4}(6e+2), 
{6e+2 \choose 5}(6e+3)
}
\right] \nonumber
\end{align}
modulo $A(12kp+10f,12e+7)+B(12kp+10f,12e+7)$ in $F(12kp+10f,12e+7)$
by Lemma~\ref{omeganoseki} (1).

We have
\begin{align}
\label{10f1-4}
0 \equiv & z_{10e,12e+1}\equiv x_{10e,12e+1}w^{6e+1} \\
\equiv & 
x_{10e,12e+1} + (6e+1)x_{10e+1,12e+2}+ {6e+1 \choose 2}x_{10e+2,12e+3}  \nonumber
\\
& + {6e+2 \choose 3}x_{10e+3,12e+4} + 
 {6e+2 \choose 4}x_{10e+4,12e+5} \nonumber \\
 &+  {6e+3 \choose 5}x_{10e+5,12e+6} +
\cdots  \nonumber \\
\equiv & \left[
{\textstyle
0, 1, 6e+1, {6e+1 \choose 2}, 
{6e+2 \choose 3}, {6e+2 \choose 4}, 
{6e+3 \choose 5}
}
\right] \nonumber
\end{align}
modulo $A(12kp+10f,12e+7)+B(12kp+10f,12e+7)$ in $F(12kp+10f,12e+7)$
by Lemma~\ref{omeganoseki} (1).

We have
\begin{align}
\label{10f1-5}
0 \equiv & z_{10e+1,12e+2} \equiv x_{10e+1,12e+2}w^{6e+1} \\
\equiv & 
x_{10e+1,12e+2} + (6e+1)x_{10e+2,12e+3}+ {6e+2 \choose 2}x_{10e+3,12e+4}  \nonumber
\\
& + {6e+2 \choose 3}x_{10e+4,12e+5} + 
 {6e+3 \choose 4}x_{10e+5,12e+6} +\cdots \nonumber \\
\equiv & \left[
{\textstyle
0, 0,1, 6e+1, {6e+2 \choose 2}, 
{6e+2 \choose 3}, {6e+3 \choose 4}
}
\right] \nonumber
\end{align}
modulo $A(12kp+10f,12e+7)+B(12kp+10f,12e+7)$ in $F(12kp+10f,12e+7)$
by Lemma~\ref{omeganoseki} (2).

We have
\begin{align}
\label{10f1-6}
0 \equiv & z_{10e+2,12e+3} \equiv x_{10e+2,12e+3}w^{6e+2} \\
\equiv & 
x_{10e+2,12e+3} + (6e+2)x_{10e+3,12e+4}+ {6e+2 \choose 2}x_{10e+4,12e+5} + {6e+3 \choose 3}x_{10e+5,12e+6} +\cdots \nonumber \\
\equiv & \left[
{\textstyle
0, 0,0,1, 6e+2, {6e+2 \choose 2}, 
{6e+3 \choose 3}
}
\right] \nonumber
\end{align}
modulo $A(12kp+10f,12e+7)+B(12kp+10f,12e+7)$ in $F(12kp+10f,12e+7)$
by Lemma~\ref{omeganoseki} (1).

We have
\begin{align}
\label{10f1-7}
0 \equiv & z_{10e+3,12e+4} \equiv x_{10e+3,12e+4}w^{6e+2} \\
\equiv & 
x_{10e+3,12e+4} + (6e+2)x_{10e+4,12e+5}+ {6e+3 \choose 2}x_{10e+5,12e+6} +\cdots \nonumber \\
\equiv & \left[
{\textstyle
0, 0,0,0,1, 6e+2, {6e+3 \choose 2}
}
\right] \nonumber
\end{align}
modulo $A(12kp+10f,12e+7)+B(12kp+10f,12e+7)$ in $F(12kp+10f,12e+7)$
by Lemma~\ref{omeganoseki} (2).

We have
\begin{align}
\label{10f1-8}
0 \equiv & z_{10e+4,12e+5} \equiv x_{10e+4,12e+5}w^{6e+3} \\
\equiv & 
x_{10e+4,12e+5} + (6e+3)x_{10e+5,12e+6}+ \cdots \nonumber \\
\equiv & \left[
{\textstyle
0, 0,0,0,0,1, 6e+3
}
\right] \nonumber
\end{align}
modulo $A(12kp+10f,12e+7)+B(12kp+10f,12e+7)$ in $F(12kp+10f,12e+7)$
by Lemma~\ref{omeganoseki} (1).

Here remark $e\equiv f \mod p$.

By (\ref{10f1-2}),  (\ref{10f1-3}),  (\ref{10f1-4}),  (\ref{10f1-5}),  (\ref{10f1-6}),  (\ref{10f1-7}),  (\ref{10f1-8}), 
we obtain 
\begin{align*}
0 \equiv & \frac{(2f+1)! (3f)!}{(5f)!}(\ref{10f1-2}) -\frac{2f+1}{6f}(\ref{10f1-3})
+(2f+1)(\ref{10f1-4})
+(3f)(\ref{10f1-5}) -\frac{1+14f+18f^2}{2} (\ref{10f1-6}) \\
& +\frac{4+53f+108f^2+63f^3}{6}(\ref{10f1-7}) -\frac{(1+12f+15f^2)(4+5f+3f^2)}{6} (\ref{10f1-8}) \\
= & \left[0,0,0,0,0,0,\frac{-24 - 308 f - 600 f^2 - 205 f^3 + 114 f^4 - 117 f^5}{40}\right] \\
= & \frac{-(10f+1)(2501757 + 5782430 f + 2175700 f^2 - 1257000 f^3 + 1170000 f^4)+101757}{4\times 10^6}x_{10e+5,12e+6}
\end{align*}
modulo $A(12kp,12e+7)+B(12kp,12e+7)$ in $F(12kp,12e+7)$.
Since $101757 = 3\times 107 \times 317$, it is not equivalent to $0$ modulo $p$.
(Recall that $p = 10f+1 \equiv 1 \mod 10$.)
Thus we obtain (\ref{10f1-0}).

For $h>0$, taking the $p^{h-1}$th power of the above equation,
we obtain
\[
x_{10ep^{h-1}+5p^{h-1},12ep^{h-1}+6p^{h-1}}\equiv 0
\]
modulo $A(12kp^h,12ep^{h-1}+6p^{h-1}+1)+B(12kp^h,12ep^{h-1}+6p^{h-1}+1)$
in $F(12kp^h,12ep^{h-1}+6p^{h-1}+1)$ by (\ref{xnoseki}).
Here remark 
\begin{align*}
10ep^{h-1}+5p^{h-1} = & 10kp^h+10fp^{h-1}+5p^{h-1} \\
= & 10\left( kp^h + fp^{h-1}+\frac{p-1}{2}(p^{h-2} + p^{h-3}+ \cdots + p+1)  \right) +5 , \\
12ep^{h-1}+6p^{h-1} = & 12kp^h+12fp^{h-1}+6p^{h-1} \\
= & 12\left( kp^h + fp^{h-1}+ \frac{p-1}{2}(p^{h-2} + p^{h-3}+ \cdots + p+1)  \right) +6 .
\end{align*}
Therefore we know
\begin{equation}\label{pnobaisuu4}
\mbox{$\maru{2}_{kp^h + fp^{h-1}+5f(p^{h-2} + p^{h-3}+ \cdots + p+1) }$ is in $D_{p^r,j}$ if $r \ge h>0$, $jp^r \le kp^h < (j+1)p^r$ and $(k,p)=1$.}
\end{equation}
Then $D_{p^r,j}$ contains $p^r-p^{r-1}$ elements of the form $\maru{1}_d$ as in (\ref{tagainiso2}) and $p^{r-1}$ elements of the form $\maru{2}_d$ as in (\ref{pnobaisuu4}).
Therefore $D_{p^r,j}$ contains at least $p^r$ elements.
We know that ${\rm Cox}(Y)$ is not Noetherian by (\ref{ets}).

\vspace{3mm}

\noindent
[\RMN{3}-4] \
Assume that $p=10f-1$, where $f >1$.

Let $k$ be a positive integer such that $(k,p)=1$.
Put $e=kp+f$.

First we shall prove that
\begin{equation}\label{10f-1-0}
x_{10e+5,12e+6} \equiv 0
\end{equation}
modulo $A(12kp,12e+7)+B(12kp,12e+7)$ in $F(12kp,12e+7)$.

Consider the triangle $T'$ with three vertices $(10e-1,12kp+p)$,
$(10e-1,12e+6)$, $(10e+6+2f,12e+6)$.
Here remark $10kp+p=10e-1$.
\begin{equation}\label{pic3}
\begin{tikzpicture}[xscale = 0.45, yscale = 0.45] 
\draw[color=gray] (0.5,-0.5) grid(15.5,13.5);
\draw[line width=0.8pt] (0.5,-0.5) grid[step=5](15.5,13.5);
%\draw[->,>=stealth,semithick] (-6,0)--(22,0)node[below]{}; %x軸
%\draw[->,>=stealth,thick] (0,-2)--(0,15)node[right]{}; %y軸
\fill (5,5) circle[radius=0.3];
\draw[color=green, line width=1pt]  (0.5,12)--(14.5,12);
\draw (0.5,12)node[left]{$12e+7$}; 
\draw (0.5,11)node[left]{$12e+6$}; 
\draw (0.5,5)node[left]{$12e$}; 
\draw (0.5,0)node[left]{$12kp+p$}; 
\draw[color=blue, line width=1pt]  (10,11) circle[radius=0.3];
\draw[color=blue, line width=1pt]  (6,6) circle[radius=0.3];
\draw (5,5) node{$\bullet$};
\foreach \x in {0,1,...,9} \draw (4+\x,2+\x) node{$\bullet$};
\foreach \x in {0,1,...,8} \draw (4+\x,3+\x) node{$\bullet$};
\foreach \x in {0,1,...,4} \draw (7+\x,7+\x) node{$\bullet$};
\foreach \x in {0,1,...,11} \draw (4+\x,\x) node{$\bullet$};
\foreach \x in {0,1,...,10} \draw (4+\x,1+\x) node{$\bullet$};
\draw[color=red] (5,11) node{$\bullet$};
\draw[color=red] (5,5) node{$\bullet$};
\foreach \x in {0,1} \draw[color=red] (5+\x,10+\x) node{$\bullet$};
\foreach \x in {0,1,2} \draw[color=red] (5+\x,9+\x) node{$\bullet$};
\foreach \x in {0,1,...,3} \draw[color=red] (5+\x,8+\x) node{$\bullet$};
\foreach \x in {0,1,...,4} \draw[color=red] (5+\x,7+\x) node{$\bullet$};
\foreach \x in {0,1,...,4} \draw[color=red] (5+\x,6+\x) node{$\bullet$};
\foreach \x in {0,1,...,7} \draw[color=red] (4,4+\x) node{$\bullet$};
\draw (5,-0.5)node[below]{$10e$}; 
\draw (15,-0.5)node[below]{$10e+6+2f$}; 
 \end{tikzpicture}
\end{equation}
For $d_1,d_2,c_0,c_1\ldots,c_6 \in K$, we put
\begin{align*}
 [d_1,d_2,c_0,c_1\ldots,c_6]
= & d_1x_{10e-1,12e-1}+d_2x_{10e+1,12e+1}+c_0x_{10e-1,12e}+c_1x_{10e,12e+1}+c_2x_{10e+1,12e+2}\\
& +c_3x_{10e+2,12e+3} +c_4x_{10e+3,12e+4}+c_5x_{10e+4,12e+5}+c_6x_{10e+5,12e+6}
\end{align*}
If $(\alpha,n)$ is in the triangle $T'$ such that $n \le \alpha + 2e$,
then $x_{\alpha,n}$ is in $A(12kp+p,12e+7)$ except for $x_{10e-1,12e-1}$ and $x_{10e+1,12e+1}$.
If $(\alpha,n)$ is in the triangle $T'$ such that $n \ge \alpha + 2e+2$,
then $x_{\alpha,n}$ is in $B(12kp+p,12e+7)$.
Therefore any $K$-linear combination of $x_{\alpha,n}$'s in the triangle $T'$ is equivalent to some  $[d_1,d_2,c_0,c_1,\ldots,c_6]$ 
modulo $A(12kp+p,12e+7)+B(12kp+p,12e+7)$ in $F(12kp+p,12e+7)$.

Taking the $p$th power of (\ref{tagainiso1}), we obtain 
\[
x_{10kp,12kp}(vx)^p\equiv 0
\]
modulo $A(12kp,12e+7)+B(12kp,12e+7)$ in $F(12kp,12e+7)$.
Here remark $12kp+2p=12kp+20f-2>12kp+12f+7=12e+7$ because $f>1$.
We know that $x_{10kp,12kp}(vx)^p \mod x^{12e+7}F$ is a $K$-linear combination of
$x_{\alpha,n}$'s in the triangle $T'$.
We have
\begin{align}
\label{10f-1-2}
& x_{10kp,12kp}(vx)^p=x_{10kp+p,12kp+p}w^{5f-1} = x_{10e-1,12e-1-2f}w^{5f-1} \\
= & 
 (1-x)^{5f-1}x_{10e-1,12e-1-2f} + (5f-1)(1-x)^{5f-1}x_{10e,12e-2f}+ {5f \choose 2}(1-x)^{5f-2}x_{10e+1,12e+1-2f}  \nonumber
\\
& + {5f \choose 3}(1-x)^{5f-2}x_{10e+2,12e+2-2f} + 
 {5f+1 \choose 4}(1-x)^{5f-3}x_{10e+3,12e+3-2f} + 
 \nonumber \\
& + {5f+1 \choose 5}(1-x)^{5f-3}x_{10e+4,12e+4-2f} + 
 {5f+2 \choose 6}(1-x)^{5f-4}x_{10e+5,12e+5-2f} +  \cdots \nonumber \\
\equiv & \left[
{\textstyle
{5f-1 \choose 2f}, {5f \choose 2}{5f-2 \choose 2f}, -{5f-1 \choose 2f+1}, -(5f-1){5f-1 \choose 2f+1}, -{5f \choose 2}{5f-2 \choose 2f+1}, 
-{5f \choose 3}{5f-2 \choose 2f+1}, 
}\right.  \nonumber \\
& \left. {\textstyle -{5f+1 \choose 4}{5f-3 \choose 2f+1}, -{5f+1 \choose 5}{5f-3 \choose 2f+1}, -{5f+2 \choose 6}{5f-4 \choose 2f+1}
}
\right] \nonumber
\end{align}
modulo $A(12kp+p,12e+7)+B(12kp+p,12e+7)$ in $F(12kp+p,12e+7)$
by Lemma~\ref{omeganoseki} (2).

We have
\begin{align}
\label{10f-1-3}
0 \equiv & z_{10e-1,12e-1}\equiv x_{10e-1,12e-1}w^{6e-1} (1+x+x^2+\cdots)^{12e-1} \\
\equiv & \left. \left( 1+(12e-1)x \vphantom{{6e \choose 2}} \right) 
\right(
(1-x)^{6e-1}x_{10e-1,12e-1} + (6e-1)(1-x)^{6e-1}x_{10e,12e} \nonumber
\\
& + {6e \choose 2}(1-x)^{6e-2}x_{10e+1,12e+1}+ {6e \choose 3}(1-x)^{6e-2}x_{10e+2,12e+2}  \nonumber \\
 &+ 
 {6e+1 \choose 4}(1-x)^{6e-3}x_{10e+3,12e+3}+  {6e+1 \choose 5}(1-x)^{6e-3}x_{10e+4,12e+4}  \nonumber \\
& \left. +
 {6e+2 \choose 6}(1-x)^{6e-4}x_{10e+5,12e+5} +
\cdots \right) \nonumber \\
\equiv & \left[
{\textstyle
1, {6e \choose 2}, 6e, (6e-1)6e, {6e \choose 2}(6e+1), 
{6e \choose 3}(6e+1), {6e+1 \choose 4}(6e+2), 
{6e+1 \choose 5}(6e+2), {6e+2 \choose 6}(6e+3)
}
\right] \nonumber
\end{align}
modulo $A(12kp+p,12e+7)+B(12kp+p,12e+7)$ in $F(12kp+p,12e+7)$
by Lemma~\ref{omeganoseki} (2).

We have
\begin{align}
\label{10f-1-4}
0 \equiv & z_{10e,12e}=
x_{10e,12e} w^{6e}(1+x+x^2+\cdots)^{12e} \\
= & (1+12ex)\left(
(1-x)^{6e}x_{10e,12e} + {6e \choose 1}(1-x)^{6e-1}x_{10e+1,12e+1}+ {6e \choose 2}(1-x)^{6e-1}x_{10e+2,12e+2} \right. \nonumber
\\
& + {6e+1 \choose 3}(1-x)^{6e-2}x_{10e+3,12e+3} + 
 {6e+1 \choose 4}(1-x)^{6e-2}x_{10e+4,12e+4} \nonumber \\
 & \left. +  {6e+2 \choose 5}(1-x)^{6e-3}x_{10e+5,12e+5} +
\cdots  \right) \nonumber \\
\equiv & \left[
{\textstyle
0,6e, 0,6e, 6e(6e+1), {6e \choose 2}(6e+1), 
{6e+1 \choose 3}(6e+2), {6e+1 \choose 4}(6e+2), 
{6e+2 \choose 5}(6e+3)
}
\right] \nonumber
\end{align}
modulo $A(12kp+p,12e+7)+B(12kp+p,12e+7)$ in $F(12kp+p,12e+7)$
by Lemma~\ref{omeganoseki} (1).

We have
\begin{align}
\label{10f-1-5}
0 \equiv & z_{10e-1,12e}\equiv
x_{10e-1,12e} w^{6e} \\
= & 
x_{10e-1,12e} + {6e \choose 1}x_{10e,12e+1}+ {6e+1 \choose 2}x_{10e+1,12e+2} + {6e+1 \choose 3}x_{10e+2,12e+3} \nonumber
\\
 &+ 
 {6e+2 \choose 4}x_{10e+3,12e+4}+  {6e+2 \choose 5}x_{10e+4,12e+5} +  {6e+3 \choose 6}x_{10e+5,12e+6} +
\cdots  \nonumber \\
\equiv & \left[
{\textstyle
0, 0,1,6e, {6e+1 \choose 2}, 
{6e+1 \choose 3}, {6e+2 \choose 4}, 
{6e+2 \choose 5},{6e+3 \choose 6}
}
\right] \nonumber
\end{align}
modulo $A(12kp+p,12e+7)+B(12kp+p,12e+7)$ in $F(12kp+p,12e+7)$
by Lemma~\ref{omeganoseki} (2).

We have
\begin{align}
\label{10f-1-6}
0 \equiv & z_{10e,12e+1}\equiv
x_{10e,12e+1} w^{6e+1} \\
= & 
x_{10e,12e+1} + {6e+1 \choose 1}x_{10e+1,12e+2}+ {6e+1 \choose 2}x_{10e+2,12e+3} + {6e+2 \choose 3}x_{10e+3,12e+4} \nonumber
\\
 &+ 
 {6e+2 \choose 4}x_{10e+4,12e+5}+  {6e+3 \choose 5}x_{10e+5,12e+6} +
\cdots  \nonumber \\
\equiv & \left[
{\textstyle
0, 0,0,1,6e+1, {6e+1 \choose 2}, 
{6e+2 \choose 3}, {6e+2 \choose 4}, 
{6e+3 \choose 5}
}
\right] \nonumber
\end{align}
modulo $A(12kp+p,12e+7)+B(12kp+p,12e+7)$ in $F(12kp+p,12e+7)$
by Lemma~\ref{omeganoseki} (1).

We have
\begin{align}
\label{10f-1-7}
0 \equiv & z_{10e+1,12e+2}\equiv
x_{10e+1,12e+2} w^{6e+1} \\
= & 
x_{10e+1,12e+2} + {6e+1 \choose 1}x_{10e+2,12e+3}+ {6e+2 \choose 2}x_{10e+3,12e+4} + {6e+2 \choose 3}x_{10e+4,12e+5} \nonumber
\\
 &+ 
 {6e+3 \choose 4}x_{10e+5,12e+6}+
\cdots  \nonumber \\
\equiv & \left[
{\textstyle
0, 0,0,0,1,6e+1, {6e+2 \choose 2}, 
{6e+2 \choose 3}, {6e+3 \choose 4}
}
\right] \nonumber
\end{align}
modulo $A(12kp+p,12e+7)+B(12kp+p,12e+7)$ in $F(12kp+p,12e+7)$
by Lemma~\ref{omeganoseki} (2).

We have
\begin{align}
\label{10f-1-8}
0 \equiv & z_{10e+2,12e+3}\equiv
x_{10e+2,12e+3} w^{6e+2} \\
= & 
x_{10e+2,12e+3} + {6e+2 \choose 1}x_{10e+3,12e+4}+ {6e+2 \choose 2}x_{10e+4,12e+5} + {6e+3 \choose 3}x_{10e+5,12e+6} +\cdots \nonumber
\\
\equiv & \left[
{\textstyle
0, 0,0,0,0,1,6e+2, {6e+2 \choose 2}, 
{6e+3 \choose 3}
}
\right] \nonumber
\end{align}
modulo $A(12kp+p,12e+7)+B(12kp+p,12e+7)$ in $F(12kp+p,12e+7)$
by Lemma~\ref{omeganoseki} (1).

We have
\begin{align}
\label{10f-1-9}
0 \equiv & z_{10e+3,12e+4}\equiv
x_{10e+3,12e+4} w^{6e+2} \\
= & 
x_{10e+3,12e+4} + {6e+2 \choose 1}x_{10e+4,12e+5}+ {6e+3 \choose 2}x_{10e+5,12e+6}  +\cdots \nonumber
\\
\equiv & \left[
{\textstyle
0, 0,0,0,0,0,1,6e+2, {6e+3 \choose 2}
}
\right] \nonumber
\end{align}
modulo $A(12kp+p,12e+7)+B(12kp+p,12e+7)$ in $F(12kp+p,12e+7)$
by Lemma~\ref{omeganoseki} (2).

We have
\begin{align}
\label{10f-1-10}
0 \equiv & z_{10e+4,12e+5}\equiv
x_{10e+4,12e+5} w^{6e+3} \\
= & 
x_{10e+4,12e+5} + {6e+3 \choose 1}x_{10e+5,12e+6} +\cdots \nonumber
\\
\equiv & \left[
{\textstyle
0, 0,0,0,0,0,0,1,6e+3
}
\right] \nonumber
\end{align}
modulo $A(12kp+p,12e+7)+B(12kp+p,12e+7)$ in $F(12kp+p,12e+7)$
by Lemma~\ref{omeganoseki} (1).

Here remark $e\equiv f \mod p$.

By (\ref{10f-1-2}),  (\ref{10f-1-3}),  (\ref{10f-1-4}),  (\ref{10f-1-5}),  (\ref{10f-1-6}),  (\ref{10f-1-7}),  (\ref{10f-1-8}),  (\ref{10f-1-9}),  (\ref{10f-1-10}), 
we obtain 
\begin{align*}
0 \equiv & \frac{(2f+1)!(3f-1!)}{(5f-1)!}(\ref{10f-1-2}) - (1+2f) (\ref{10f-1-3})
+\frac{(1+2f)(-1+21f)}{12}(\ref{10f-1-4}) - (1-9f-12f^2)(\ref{10f-1-5}) \\
& -\frac{-2+15f+49f^2+42f^3}{2}(\ref{10f-1-6})
-\frac{2-14f-51f^2-45f^3}{2}(\ref{10f-1-7}) \\
& -\frac{-12+79f+335f^2+432f^3+198f^4}{12}(\ref{10f-1-8})
+\frac{(1+f)^2(-8+66f+177f^2+81f^3)}{8}(\ref{10f-1-9}) \\
& -\frac{-120+718f+3821f^2+6630f^3+3480f^4-5508f^5-7101f^6}{120}(\ref{10f-1-10}) \\
= & \left[0,0,0,0,0,0,0,0,
\frac{3(2 + 3f) (120 - 872f - 2670f^2 - 3213f^3 + 4907f^4 + 
   22509f^5 + 24147f^6)}{720}
\right] 
\end{align*}
modulo $A(12kp,12e+7)+B(12kp,12e+7)$ in $F(12kp,12e+7)$.
We have
\begin{align}
& 10^7\cdot 3 (2 + 3f) (120 - 872f - 2670f^2 - 3213f^3 + 4907f^4 + 
   22509f^5 + 24147f^6) \label{remainder1} \\
= & 3(20+3\cdot 10f)(12\cdot 10^7 -872\cdot  10^5(10f) -267\cdot  10^5(10f)^2
-3213\cdot  10^3(10 f)^3 \nonumber \\
& + 4907\cdot  10^2(10f)^4 + 22509\cdot  10(10f)^5 + 24147\cdot (10f)^6)  \nonumber \\
\equiv & 3 \cdot 23 \cdot 3626937=3^5\cdot 23\cdot 44777 \mod (10f-1)\nonumber .
\end{align}
Therefore it is not equivalent to $0$ modulo $p$.
(Recall that $p = 10f-1 \equiv 9 \mod 10$.)
Thus we obtain (\ref{10f-1-0}).

For $h>0$, taking the $p^{h-1}$th power of (\ref{10f-1-0}),
we obtain
\[
x_{10ep^{h-1}+5p^{h-1},12ep^{h-1}+6p^{h-1}}\equiv 0
\]
modulo $A(12kp^h,12ep^{h-1}+6p^{h-1}+1)+B(12kp^h,12ep^{h-1}+6p^{h-1}+1)$
in $F(12kp^h,12ep^{h-1}+6p^{h-1}+1)$ by (\ref{xnoseki}).
Here remark 
\begin{align*}
10ep^{h-1}+5p^{h-1} = & 10kp^h+10fp^{h-1}+5p^{h-1} \\
= & 10\left( kp^h + fp^{h-1}+\frac{p-1}{2}(p^{h-2} + p^{h-3}+ \cdots + p+1)  \right) +5 , \\
12ep^{h-1}+6p^{h-1} = & 12kp^h+12fp^{h-1}+6p^{h-1} \\
= & 12\left( kp^h + fp^{h-1}+ \frac{p-1}{2}(p^{h-2} + p^{h-3}+ \cdots + p+1)  \right) +6 .
\end{align*}
Therefore we know
\begin{equation}\label{pnobaisuu5}
\mbox{$\maru{2}_{kp^h + fp^{h-1}+(5f-1)(p^{h-2} + p^{h-3}+ \cdots + p+1) }$ is in $D_{p^r,j}$ if $r \ge h>0$, $jp^r \le kp^h < (j+1)p^r$ and $(k,p)=1$.}
\end{equation}
Then $D_{p^r,j}$ contains $p^r-p^{r-1}$ elements of the form $\maru{1}_d$ as in (\ref{tagainiso2}) and $p^{r-1}$ elements of the form $\maru{2}_d$ as in (\ref{pnobaisuu5}).
Therefore $D_{p^r,j}$ contains at least $p^r$ elements.
We know that ${\rm Cox}(Y)$ is not Noetherian by (\ref{ets}).

\vspace{3mm}

\noindent
[\RMN{3}-5] \
Assume that $p=10f+7$, where $f \ge 0$.

Let $k$ be a positive integer such that $(k,p)=1$.
Put $e=kp+f$.

First we shall prove 
\begin{equation}\label{10f7-0}
x_{10e+7,12e+8} \equiv 0
\end{equation}
modulo $A(12kp,12e+9)+B(12kp,12e+9)$ in $F(12kp,12e+9)$.
Taking the $p$th power of (\ref{tagainiso1}), we obtain 
\begin{equation}\label{10f7-1}
x_{10kp,12kp}(vx)^p \equiv 0
\end{equation}
modulo $A(12kp,12e+9)+B(12kp,12e+9)$ in $F(12kp,12e+9)$.
(Remark $12kp+2p=12kp+20f+14>12kp+12f+9=12e+9$.)
Then we have
\begin{align*}
& x_{10kp,12kp}(vx)^p = x_{10e+7,12e+7-2f}w^{5f+3} \\
= & (1-x)^{5f+3}x_{10e+7,12e+7-2f} + {5f+3 \choose 1}(1-x)^{5f+3}x_{10e+8,12e+8-2f}
+ \cdots
\end{align*}
by Lemma~\ref{omeganoseki} (2).
The coefficient of $x_{10e+7,12e+8}$ is $-{5f+3 \choose 2f+1}$.
Thus (\ref{10f7-0}) follows from this.
Therefore we know
\begin{equation}\label{pnobaisuu6}
\mbox{$\maru{3}_{kp+f }$ is in $D_{p^r,j}$ if $(k,p)=1$ and $jp^r \le kp < (j+1)p^r$.}
\end{equation}

Let $m$ be the integer satisfying 
\begin{equation}\label{44777-2}
(10e+7)p+1=10m, 
\end{equation}
that is,
\[
m=ep+7f+5=kp^2+fp+7f+5 .
\]

Next we shall prove  
\begin{itemize}
\item[(1)]
if $p \neq 44777$, then
\begin{equation}\label{10f7-2}
x_{10m+5,12m+6} \equiv 0
\end{equation}
modulo $A(12kp^2,12m+7)+B(12kp^2,12m+7)$ in $F(12kp^2,12m+7)$,
\item[(2)]
if $p= 44777$, then
\begin{equation}\label{10f7-3}
x_{10m+7,12m+8} \equiv 0
\end{equation}
modulo $A(12kp^2,12m+9)+B(12kp^2,12m+9)$ in $F(12kp^2,12m+9)$.
\end{itemize}

Consider the triangle $T''$ with three vertices $(10m-1,12ep+8p)$,
$(10m-1,12m+8)$, $(10m+11+4f,12m+8)$.
\begin{equation}\label{pic4}
\begin{tikzpicture}[xscale = 0.45, yscale = 0.45] 
\draw[color=gray] (0.5,-0.5) grid(16.5,14.5);
\draw[line width=0.8pt] (0.5,-0.5) grid[step=5](16.5,14.5);
%\draw[->,>=stealth,semithick] (-6,0)--(22,0)node[below]{}; %x軸
%\draw[->,>=stealth,thick] (0,-2)--(0,15)node[right]{}; %y軸
\fill (5,5) circle[radius=0.3];
\draw[color=green, line width=1pt]  (0.5,14)--(16.5,14);
\draw (0.5,13)node[left]{$12m+8$}; 
\draw (0.5,11)node[left]{$12m+6$}; 
\draw (0.5,5)node[left]{$12m$}; 
\draw (0.5,1)node[left]{$12ep+8p$}; 
\draw[color=blue, line width=1pt]  (10,11) circle[radius=0.3];
\draw[color=blue, line width=1pt]  (6,6) circle[radius=0.3];
\draw[color=blue, line width=1pt]  (12,13) circle[radius=0.3];
\draw (5,5) node{$\bullet$};
\draw (11,12) node{$\bullet$};
\foreach \x in {0,1,...,11} \draw (4+\x,2+\x) node{$\bullet$};
\foreach \x in {0,1,...,10} \draw (4+\x,3+\x) node{$\bullet$};
\foreach \x in {0,1,...,6} \draw (7+\x,7+\x) node{$\bullet$};
%\foreach \x in {0,1,...,11} \draw (4+\x,\x) node{$\bullet$};
\foreach \x in {0,1,...,12} \draw (4+\x,1+\x) node{$\bullet$};
\draw[color=red] (5,11) node{$\bullet$};
\draw[color=red] (5,5) node{$\bullet$};
\foreach \x in {0,1} \draw[color=red] (5+\x,12) node{$\bullet$};
\foreach \x in {0,1,2} \draw[color=red] (5+\x,13) node{$\bullet$};
\foreach \x in {0,1,2,3} \draw[color=red] (5+\x,10+\x) node{$\bullet$};
\foreach \x in {0,1,2,3,4} \draw[color=red] (5+\x,9+\x) node{$\bullet$};
\foreach \x in {0,1,...,5} \draw[color=red] (5+\x,8+\x) node{$\bullet$};
\foreach \x in {0,1,...,6} \draw[color=red] (5+\x,7+\x) node{$\bullet$};
\foreach \x in {0,1,...,4} \draw[color=red] (5+\x,6+\x) node{$\bullet$};
\foreach \x in {0,1,...,9} \draw[color=red] (4,4+\x) node{$\bullet$};
\draw (5,-0.5)node[below]{$10m$}; 
\draw (16,-0.5)node[below]{$10m+11+4f$}; 
\draw[line width=1pt] (16,-0.5)--(16,0);
 \end{tikzpicture}
\end{equation}

For $d_1,d_2,c_0,c_1\ldots,c_7 \in K$, we put
\begin{align*}
&  [d_1,d_2,c_0,c_1\ldots,c_7]\\
= & d_1x_{10m-1,12m-1}+d_2x_{10m+1,12m+1}+c_0x_{10m-1,12m}+c_1x_{10m,12m+1}+c_2x_{10m+1,12m+2} \\
& +c_3x_{10m+2,12m+3}+c_4x_{10m+3,12m+4}+c_5x_{10m+4,12m+5}+c_6x_{10m+5,12m+6}+c_7x_{10m+7,12m+8}
\end{align*}
If $(\alpha,n)$ is in the triangle $T''$ such that $n \le \alpha + 2m$,
then $x_{\alpha,n}$ is in $A(12ep+8p,12m+9)$ except for $x_{10m-1,12m-1}$ and $x_{10m+1,12m+1}$.
If $(\alpha,n)$ is in the triangle $T''$ such that $n \ge \alpha + 2m+2$,
then $x_{\alpha,n}$ is in $B(12ep+8p,12m+9)$.
Therefore any $K$-linear combination of $x_{\alpha,n}$'s in the triangle $T''$ is equivalent to some  $[d_1,d_2,c_0,c_1,\ldots,c_7]$ 
modulo $A(12ep+8p,12m+9)+B(12ep+8p,12m+9)$ in $F(12ep+8p,12m+9)$.

Since $F(12e+9,12e+10)=A(12e+9,12e+10)+B(12e+9,12e+10)$, we btain
\begin{equation}\label{10f7-0'}
x_{10e+7,12e+8} \equiv 0
\end{equation}
modulo $A(12kp^2,12e+10)+B(12kp^2,12e+10)$ in $F(12kp^2,12e+10)$
by (\ref{10f7-0}).
Taking the $p$th power of (\ref{10f7-0'}), we obtain 
\begin{align}\label{10f7-4}
0\equiv & (x_{10e+7,12e+8} )^p
=(v^{10e+7}w^{\lceil -\frac{10e+7}{2}  \rceil}x^{12e+8})^p
=x_{10m-1,12ep+8p}w^{5f+3} \\
= & 
 (1-x)^{5f+3}x_{10m-1,12ep+8p} + (5f+3)(1-x)^{5f+3}x_{10m,12ep+8p+1}  \nonumber
\\
& + {5f+4 \choose 2}(1-x)^{5f+2}x_{10m+1,12ep+8p+2}+ {5f+4 \choose 3}(1-x)^{5f+2}x_{10m+2,12ep+8p+3}  
 \nonumber \\
& + 
 {5f+5 \choose 4}(1-x)^{5f+1}x_{10m+3,12ep+8p+4}+ {5f+5 \choose 5}(1-x)^{5f+1}x_{10m+4,12ep+8p+5}   \nonumber \\
& + 
 {5f+6 \choose 6}(1-x)^{5f}x_{10m+5,12ep+8p+6} + {5f+7 \choose 8}(1-x)^{5f-1}x_{10m+7,12ep+8p+8}+ \cdots \nonumber \\
\equiv & \left[
{\textstyle
-{5f+3 \choose 4f+3}, -{5f+4 \choose 2}{5f+2 \choose 4f+3}, {5f+3 \choose 4f+4}, (5f+3){5f+3 \choose 4f+4}, {5f+4 \choose 2}{5f+2 \choose 4f+4}, 
{5f+4 \choose 3}{5f+2 \choose 4f+4}, 
}\right.  \nonumber \\
& \left. {\textstyle {5f+5 \choose 4}{5f+1 \choose 4f+4}, {5f+5 \choose 5}{5f+1 \choose 4f+4}, {5f+6 \choose 6}{5f \choose 4f+4},{5f+7 \choose 8}{5f-1 \choose 4f+4}
}
\right] \nonumber
\end{align}
modulo $A(12ep+8p,12m+9)+B(12ep+8p,12m+9)$ in $F(12ep+8p,12m+9)$
since $12ep+10p=12ep+100f+70\ge 12ep+84f+69=12m+9$.
Here remark that (\ref{10f7-4}) is $[-1,0,\cdots,0]$ if $f=0$.

Replacing $e$ by $m$ and adding the last component\footnote{If $p\neq 44777$, then we do not need the last component.
If $p = 44777$, then we may assume that denominators of ${5f+7 \choose 8}$,
${6m+3 \choose 8}$, ${6m+3 \choose 7}$, ${6m+4 \choose 8}$, ${6m+4 \choose 7}$ (in the last components) are units.
} to (\ref{10f-1-3}), (\ref{10f-1-4}), (\ref{10f-1-5}), (\ref{10f-1-6}), (\ref{10f-1-7}), (\ref{10f-1-8}), (\ref{10f-1-9}), (\ref{10f-1-10}), 
we obtain
\begin{align}
0 \equiv & \left[
{\textstyle
1, {6m \choose 2}, 6m, (6m-1)6m, {6m \choose 2}(6m+1), 
{6m \choose 3}(6m+1), {6m+1 \choose 4}(6m+2), 
{6m+1 \choose 5}(6m+2),} \right.  \label{10f7-5}\\
& \left. {\textstyle
{6m+2 \choose 6}(6m+3), {6m+3 \choose 8}(6m+4)
}
\right] \nonumber\\
0 \equiv & \left[
{\textstyle
0,6m,0,6m, 6m(6m+1), {6m \choose 2}(6m+1), 
{6m+1 \choose 3}(6m+2), {6m+1 \choose 4}(6m+2), 
{6m+2 \choose 5}(6m+3),} \right.  \label{10f7-6}\\
& \left. {\textstyle
{6m+3 \choose 7}(6m+4)
}
\right] \nonumber \\
0 \equiv & \left[
{\textstyle
0, 0,1,6m, {6m+1 \choose 2}, 
{6m+1 \choose 3}, {6m+2 \choose 4}, 
{6m+2 \choose 5},{6m+3 \choose 6}, {6m+4 \choose 8}
}
\right]  \label{10f7-7}\\
0 \equiv & \left[
{\textstyle
0, 0,0,1,6m+1, {6m+1 \choose 2}, 
{6m+2 \choose 3}, {6m+2 \choose 4}, 
{6m+3 \choose 5},{6m+4 \choose 7}
}
\right] \label{10f7-8} \\
0 \equiv & \left[
{\textstyle
0, 0,0,0,1,6m+1, {6m+2 \choose 2}, 
{6m+2 \choose 3}, {6m+3 \choose 4}, 
{6m+4 \choose 6}
}
\right] \label{10f7-9} \\
0 \equiv & \left[
{\textstyle
0, 0,0,0,0,1,6m+2, {6m+2 \choose 2}, 
{6m+3 \choose 3}, {6m+4 \choose 5}
}
\right] \label{10f7-10} \\
0 \equiv & \left[
{\textstyle
0, 0,0,0,0,0,1,6m+2, {6m+3 \choose 2}, 
{6m+4 \choose 4}
}
\right] \label{10f7-11} \\
0 \equiv & \left[
{\textstyle
0, 0,0,0,0,0,0,1,6m+3, {6m+4 \choose 3}
}
\right] \label{10f7-12}
\end{align}
modulo $A(12kp+p,12m+7)+B(12kp+p,12m+7)$ in $F(12kp+p,12m+7)$.
Remark that 
\begin{equation}\label{44777-4}
m\equiv 7f+5\mod p .
\end{equation}
By (\ref{10f7-4}), (\ref{10f7-5}), (\ref{10f7-6}), (\ref{10f7-7}), (\ref{10f7-8}), (\ref{10f7-9}), (\ref{10f7-10}), (\ref{10f7-11}), (\ref{10f7-12}), 
we have
\begin{align*}
0 \equiv & \frac{(4f+4)!f!}{(5f+3)!}(\ref{10f7-4}) +(4f+4)(\ref{10f7-5})
-\frac{2(1+f)(870+2474f+1759f^2)}{6m}(\ref{10f7-6}) \\
&-(120+289f+168f^2)(\ref{10f7-7})
 +(1860+7003f+8671f^2+3518f^3)(\ref{10f7-8}) \\
& -\frac{3720+14708f+19303f^2+8405f^3}{2}(\ref{10f7-9}) \\
& +\frac{65100+343772f+680605f^2+598882f^3+197669f^4}{6}(\ref{10f7-10}) \\
& -\frac{476160+3042856f+7850330f^2+10243855f^3+6774982f^4+1819801f^5}{24}(\ref{10f7-11}) \\
& - {\textstyle \frac{190940160+1516966184f+5005402706f^2+8777761165f^3+8625823355f^4+
4502014011f^5+974544899f^6}{120} }(\ref{10f7-12})
\\
=& [0,0,0,0,0,0,0,0,q_1,q_2]
\end{align*}
modulo $A(12ep+8p,12m+9)+B(12ep+8p,12m+9)$ in $F(12ep+8p,12m+9)$, where
\begin{align}
720q_1 = & -26767572480 - 246120200736 f - 967942897272 f^2 - 
 2110412205706 f^3 \nonumber \\
 & - 2754630615405 f^4 - 2152135097539 f^5 - 
 931716713643 f^6 - 172390143619 f^7 \nonumber \\
40320q_2  = & -348081961328640 - 4085017940012352 f - 21279829406091360f^2  
\label{44777-1} \\
& - 64577996264481356 f^3
 - 125811467012647820 f^4 - 
 163172345721567295 f^5 \nonumber \\ 
 & - 140876419259495720 f^6 -  78068028418279174f^7
 - 25195471807991660 f^8 \nonumber \\ 
 &- 3607880835288623 f^9 .\nonumber 
\end{align}
Here $720$ is not divided by $p=10f+7$, 
and $40320$ is not divided by $44777$.

First assume that $p\neq 44777$.
When we divide $720\times 10^7q_1$ by $10f+7$ (as a polynomial of $f$), the reminder is $-250258653$.
Here we have 
\begin{equation}\label{remainder2}
250258653=3^5\times 23\times 44777 .
\end{equation}
Therefore (\ref{10f7-2}) holds.
(We have to give an attention to  the case $f=0$, since some binomial coefficients are $0$.)

Taking the $p^{h-2}$th power of  (\ref{10f7-2}) for $h \ge 2$,
we obtain
\[
x_{10mp^{h-2}+5p^{h-2},12mp^{h-2}+6p^{h-2}}\equiv 0
\]
modulo $A(12kp^h,12mp^{h-2}+6p^{h-2}+1)+B(12kp^h,12mp^{h-2}+6p^{h-2}+1)$
in $F(12kp^h,12mp^{h-2}+6p^{h-2}+1)$ by (\ref{xnoseki}).
Here remark 
\begin{align*}
10mp^{h-2}+5p^{h-2} = & 10(kp^2+fp+7f+5)p^{h-2}+5(p^{h-2}-1)+5 \\
= & 10\left( kp^h + fp^{h-1}+(7f+5)p^{h-2}+\frac{p-1}{2}(p^{h-3} + p^{h-4}+ \cdots + p+1)  \right) +5 , \\
12mp^{h-2}+6p^{h-2} = & 12(kp^2+fp+7f+5)p^{h-2}+6(p^{h-2}-1)+6 \\
= & 12\left( kp^h + fp^{h-1}+(7f+5)p^{h-2}+\frac{p-1}{2}(p^{h-3} + p^{h-4}+ \cdots + p+1)  \right) +6  .
\end{align*}
Therefore we know
\begin{equation}\label{pnobaisuu7}
\mbox{$\maru{2}_{kp^h + fp^{h-1}+(7f+5)p^{h-2} +(5f+3)(p^{h-3} + p^{h-4}+ \cdots + p+1) }$ is in $D_{p^r,j}$ if $r \ge h\ge 2$ and $jp^r \le kp^h < (j+1)p^r$.}
\end{equation}
Then $D_{p^r,j}$ contains $p^r-p^{r-1}$ elements of the form $\maru{1}_d$ as in (\ref{tagainiso2}), $p^{r-1}-p^{r-2}$ elements of the form $\maru{3}_d$ as in (\ref{pnobaisuu6}) and
$p^{r-2}$ elements of the form $\maru{2}_d$ as in (\ref{pnobaisuu7}).
Therefore $D_{p^r,j}$ contains at least $p^r$ elements.
We know that ${\rm Cox}(Y)$ is not Noetherian by (\ref{ets}).

Finally assume $p=44777$.
When we divide $40320\times 10^9q_2$ by $10f+7$ (as a polynomial of $f$), the reminder is $-5257057765239$.
Here $5257057765239$ is not divided by $44777$.
Therefore (\ref{10f7-3}) holds.
Recall that (\ref{10f7-0}) implies (\ref{10f7-3}) for $m$ satisfying (\ref{44777-2}).
Therefore, putting
\[
(10m+7)p+1=10n, 
\]
(\ref{10f7-3}) implies 
\[
x_{10n+7,12n+8} \equiv 0
\]
modulo $A(12kp^3,12n+9)+B(12kp^3,12n+9)$ in $F(12kp^3,12n+9)$.
Here
\[
n = mp+7f+5=kp^3+fp^2+(7f+5)p+(7f+5) .
\]
Repeating this process, we know
\begin{equation}\label{pnobaisuu8}
\mbox{$\maru{3}_{kp^h + fp^{h-1}+(7f+5)(p^{h-2} + p^{h-3} + p^{h-4}+ \cdots + p+1) }$ is in $D_{p^r,j}$ if $r \ge h\ge 1$ and $jp^r \le kp^h < (j+1)p^r$.}
\end{equation}
Then $D_{p^r,j}$ contains $p^r-p^{r-1}$ elements of the form $\maru{1}_d$ as in (\ref{tagainiso2}), $p^{r-1}$ elements of the form $\maru{3}_d$ as in (\ref{pnobaisuu8}).
Therefore $D_{p^r,j}$ contains at least $p^r$ elements.
We know that ${\rm Cox}(Y)$ is not Noetherian by (\ref{ets}).

\begin{Remark}
\begin{rm}
We put $f_1(x)=2x+1$, $f_2(x)=\frac{5x}{2}(2x+1)(3x-1)$,
$f_3(x)=-(3x-1)$, $f_4(x)=-(5x-1)(3x-1)$, $f_5(x)=-\frac{5x}{2}(3x-1)(3x-2)$.
$f_6(x)=-\frac{5x}{6}(5x-2)(3x-1)(3x-2)$, $f_7(x)=-\frac{(5x+1)}{24}(5x)(3x-1)(3x-2)(3x-3)$,
$f_8(x)=-\frac{(5x+1)}{120}(5x)(5x-3)(3x-1)(3x-2)(3x-3)$,
$f_9(x)=-\frac{(5x+2)}{720}(5x+1)(5x)(3x-1)(3x-2)(3x-3)(3x-4)$,
$f_{10}(x)=-\frac{(5x+3)}{8!}(5x+2)(5x+1)(5x)(3x-1)(3x-2)(3x-3)(3x-4)(3x-5)$.

Then we have
\begin{align*}
\frac{(2f+1)!(3f-1!)}{(5f-1)!}(\ref{10f-1-2})& =
[f_1(e),f_2(e),f_3(e),f_4(e),f_5(e),f_6(e),f_7(e),f_8(e),f_9(e)], \\
-\frac{(4f+4)!f!}{(5f+3)!}(\ref{10f7-4}) & =
[f_1(m),f_2(m),f_3(m),f_4(m),f_5(m),f_6(m),f_7(m),f_8(m),f_9(m),f_{10}(m)].
\end{align*}
Furthermore we have
\begin{align*}
10e \equiv 1 \mod p(=10f-1), \\
10m\equiv 1 \mod p(=10f+7).
\end{align*}
Hence we obtained the same remainder in (\ref{remainder1}) and (\ref{remainder2}).
\end{rm}
\end{Remark}

\vspace{3mm}

\noindent
\begin{tabular}{l}
Kazuhiko Kurano \\
Department of Mathematics \\
Faculty of Science and Technology \\
Meiji University \\
Higashimita 1-1-1, Tama-ku \\
Kawasaki 214-8571, Japan \\
{\tt kurano@meiji.ac.jp} \\
{\tt http://www.isc.meiji.ac.jp/\~{}kurano}
\end{tabular}

\end{document}